\documentclass{amsart}
\usepackage{amssymb}
\usepackage{amsfonts}
\usepackage{amstext}
\usepackage{algorithmic}
\usepackage{algorithm}
\usepackage{graphicx}
\usepackage{epstopdf}
\usepackage[all]{xy}

\numberwithin{equation}{section}
\newtheorem{theorem}{Theorem}[section]
\newtheorem{lemma}[theorem]{Lemma}
\newtheorem{proposition}[theorem]{Proposition}
\newtheorem{corollary}[theorem]{Corollary}

\theoremstyle{definition}
\newtheorem{definition}[theorem]{Definition}
\newtheorem{example}[theorem]{Example}

\theoremstyle{remark}

\renewcommand{\span}{{\rm span}}

\newcommand{\R}{{\mathbb{R}}}

\newcommand{\C}{{\mathbb{C}}}
\newcommand{\Z}{{\mathbb{Z}}}

\newcommand{\CC}{{\mathcal{C}}}

\newcommand{\CL}{{\mathcal{L}}}

\renewcommand{\(}{{(}}
\renewcommand{\)}{{)}}
\newcommand{\cg}{{\mathfrak{g}}}

\renewcommand{\ker}{{\rm{ker}}}

\newcommand{\<}{{\langle}}
\renewcommand{\>}{{\rangle}}

\newcommand{\tens}{\otimes}
\newcommand{\id}{{\rm id}}
\newcommand{\extd}{{\rm d}}

\newcommand{\from}{{\leftarrow}}

\begin{document}

\title{Noncommutative Riemannian geometry on graphs}
\keywords{graph theory, Laplace-Beltrami, noncommutative geometry, Yang-Baxter, Cayley graph, finite group. Revised Dec 2011}


\author{Shahn Majid}
\address{Queen Mary University of London\\
School of Mathematics, Mile End Rd, London E1 4NS, UK}

\email{s.majid@qmul.ac.uk}


\begin{abstract} We show that arising out of noncmmutatve geometry is a natural family of {\em edge Laplacians} on the edges of a graph. The family includes a canonical edge Laplacian associated to the graph, extending the usual graph Laplacian on vertices, and we find its spectrum. We show that for a connected graph its eigenvalues are strictly positive aside from one mandatory zero mode, and include all the vertex degrees. Our edge Laplacian is not the graph Laplacian on the line graph but rather it arises as the noncommutative Laplace-Beltrami operator on differential 1-forms, where we use the language of differential algebras to functorially interpret a graph as providing a  `finite manifold structure' on the set of vertices. We equip any graph with a canonical  `Euclidean metric' and a canonical bimodule connection, and in the case of a Cayley graph we construct a metric compatible connection for the Euclidean metric.  We make use of results on bimodule connections on inner calculi on algebras, which we prove, including a general  relation between zero curvature and the braid relations. 
\end{abstract}
\maketitle

\section{Introduction}

A differential algebra is an algebra $A$ equipped with a pair $(\Omega^1,\extd)$ where $\Omega^1$ is an $A-A$-bimodule and $\extd:A\to \Omega^1$ is a linear map obeying the Leibniz rule $\extd(ab)=(\extd a)b+a\extd b$ for all $a,b\in A$. It is usually required that the map $A\tens A\to \Omega^1$ given by $a\tens b\mapsto a\extd b$ is surjective. This notion features in essentially all approaches to noncommutative geometry and has been applied extensively in the case where $A$ is noncommutative, such as to a Lie theory of quantum groups.  The space $\Omega^1$ plays the role of differentials or 1-forms in differential geometry, but because we do not suppose that the left and right module structures on $A$ are equal (i.e. 1-forms may not commute with functions)  this notion is fundamentally more general than conventional differential geometry when specialised to commutatve algebras. In particular, it  is exactly what is needed to provide a notion of differential geometry on finite sets, the only ordinary differentiable structure on a discrete topology being the zero one.

Our first result, Theorem~3.1, is to make this more precise. If $X$ is a finite set it is known that differential structures on $A=k(X)$, $k$ a field, are in 1-1 correspondence with digraphs with vertex set $X$. This is known, see for example\cite{Connes}, but we show that this correspondence is functorial. This means that natural constructions for digraphs can be expressed in terms of differential algebra and vice versa differential algebra constructions can be specialised. Although many constructions in discrete mathematics are loosely motivated by geometric intuition the precise nature of our correspondence allows one to systematically transfer those ides of classical differential geometry that can be extended to general differential algebras and then specialised. To this end much of differential geometry, notably metrics $\cg\in \Omega^1\tens_A\Omega^1$ and linear connections $\nabla:\Omega^1\to \Omega^1\tens_A\Omega^1$ on bimodules, their curvature $R_\nabla$ and geometric torsion $T_\nabla$ are all understood at the level of differential algebra an can be specialised to digraphs, which we do. This framework material is used in the noncommutative geometry literature and we give a short account in Sections~2.1-2.2, with some further algebraic results in the rest of Section~2 that will be needed later. Notably, Theorem~2.1 relates zero curvature to the braid relations for a `generalised flip' map $\sigma$ defined by a connection. 

Focussing mainly on the symmetric or `bidirected' case (i.e. we view an undirected graph as directed both ways) Section~3.2 analyses the most general form of metric and bimodule connection. Notably in Proposition~3.7 we reproduce the general graph Laplacian on vertices as $(\ ,\ )\nabla\extd$ in terms of an inverse metric  and a bimodule connection, providing a more extended geometric picture than hitherto available. Section 3.3 provides a canonical $\cg,\nabla$ which are nevertheless quite unusual from the point of view of classical geometry in that the `flip map' $\sigma$ associated to the connection is the identity. Applying the preceding analysis, Proposition~3.10 finds a natural zero-curvature `Maurer-Cartan connection' $\nabla$ on any Cayley graph as a member of a class of connections of `permutation type' studied in Section~3.5, including their torsion $T_\nabla$ and curvature $R_\nabla$.

We then come in Section~4 to the main result of the paper, a Laplacian on the edges of a graph. In classical Riemannian geometry the Hodge-de Rham theory provides for a Laplace-Beltrami (and a Hodge-Laplace) operator on all degrees of forms, not just functions, and in Section~4.1 we develop this in degree 1 at  the level of inner differential algebras and connections $\nabla$. We then  see in Section 4.2 how a natural Laplacian emerges on the edges of the graph as the analogue of the Laplace-Beltrami on the space of 1-forms on a manifold.   Section~4.3 specialises to the canonical connection and Euclidean metric of Section~3.3 and we study this `canonical edge Laplacian' in detail. We prove (Theorem~4.9) that the eigenvalues of the canonical edge Laplacian are strictly positive  for any connected graph, aside from a single zero mode (just as for the usual graph Laplacian on the vertices and extending that). We show that its spectrum has one part which is (twice) that of the usual graph Laplacian on vertices and a second part consisting of the integer degrees of every vertex. Our analysis does not exactly tell us that the edge Laplacian is diagonalisable but it typically is, a sufficient condition we prove being that the two parts of the spectrum are disjoint.

Also, going the other way, digraph geometry provides a good illustration of non-classical ideas in noncommutative geometry, being fundamentally noncommutative even though the algebra $k(V)$ is commutative. The most important of these is that sufficiently noncommutative geometries tend to be `inner' in the sense of a differential 1-form $\theta\in\Omega^1$ such that $[\theta, a]=\extd a$ for all $a\in A$ (here $[\theta,a]=\theta a-a\theta$). There can be no such concept in classical differential geometry as 1-forms and functions commute. But if classical geometry is a limit of a noncommutative geometry then there can be phenomena in classical geometry which are unconnected but which become connected or explicable in the noncommutative case. We illustrate this with the Laplacian which we show (Proposition~4.3) has a deeper origin as the `partial derivative' conjugate to the direction of $\theta$. 

\subsubsection*{Acknowledgemts} I would like to thank Carsten Thomassen for discussions which improved the paper, notably the proof of part (3) of Proposition~3.2.

\section{General framework}

We briefly describe a formulation of non-commutative  Riemannian differential algebras that we will applied in the rest of the paper.  Sections~2.1 and~2.2 are not intended to contain anything new. The other sections contain some new material for Laplacians and geometry with inner calculi, motivated by  \cite{Ma:bh}. 

\subsection{}  Let $k$ be our ground field of characteristic not 2 and  let $A$ be a unital algebra over $k$ viewed as `coordinate algebra' or functions on a space (it can, however, be noncommutative).  Throughout the paper, $\bar{\tens}=\tens_A$ for brevity. A differential algebra structure on $A$ (or `differential structure') means a specification of a space $\Omega^1$ of `1-forms'  and a linear map $\extd$ or `exterior derivative' obeying
\begin{enumerate} 
\item $\Omega^1$ is an $A-A$- bimodule   (so $f(\omega g=(f\omega)g$ for all $f,g \in A,\ \omega\in\Omega^1$
\item $\extd: A\to \Omega^1$ is a derivation $\extd(fg)=(\extd f)g+f\extd g$
\item $A\tens A\to\Omega^1$ by $f\tens g\mapsto f\extd g$ is surjective
\item (optional connectedness condition) $\ker\extd=k 1$.
\end{enumerate}

A morphism of differential algebras $(A,\Omega^1,\extd)\to (B,\Omega^1,\extd)$ means an algebra map $\phi:A\to B$ and a compatible $\phi_*:\Omega^1(A)\to \Omega^1(B)$ such that
\begin{equation}\label{moralg}  \begin{array}{rcl} \Omega^1(A) &\  \ {\buildrel \phi_* \over \longrightarrow} \ \ & \Omega^1(B)\\ \extd\nwarrow & & \nearrow\extd\\ A &{\buildrel \phi\over \longrightarrow } & B\end{array}\end{equation}
commutes. The surjectivity axiom (3) means that $\phi_*$ if it exists is uniquely determined once $\phi$ is specified, i.e. this diagram says what it means for a map between differentiable algebras to be differentiable. 

We say that a calculus is left(right) parallelizable if it is free as a left (right) $A$-module. In this case the smallest cardinality of a basis over $A$ is called the left (right) cotangent dimension. At least for $A$ trace class or commutative, all left (right) bases if they exist have the same cardinality. 

Every algebra $A$ has a connected `universal differential calculus' where $\Omega^1\subset A\tens A$ is the kernel of the product map and $\extd a=1\tens a-a\tens 1$. Every other differential structure on $A$ is a quotient of the universal one by a sub-bimodule. If $A$ is finite dimensional then the universal differential calculus is left and right parallelizable with cotangent dimension $\dim(A)-1$. 

Finally, it is always possible, not necessarily uniquely, to extend $\Omega^1$ to a differential graded algebra or `exterior algebra' $\Omega=\oplus_n\Omega^n$ where $\Omega^0=A$, $\Omega$ is generated by degree $0,1$ and $\extd$ extends as a graded derivation with $\extd^2=0$. The top degree if there is one is called the `volume dimension'. The volume dimension and the cotangent dimension are not usually the same for a general algebra, even a commutative one. The cohomology of this complex is
called the noncommutative de Rham cohomology $H(A, \Omega,\extd)$. Here $H^0=k1$ expresses that the calculus is connected. 

In particular there is a maximal prolongation of any $(A,\Omega^1,\extd)$ built on the tensor algebra of $\Omega^1$ over $A$  modulo relations required by the properties of $\extd$.   This implies that $H^1(A,\Omega,\extd)$ with the maximal prolongation is an invariant of first order differential algebras. The maximal prolongation of the universal calculus is the universal differential graded algebra on $A$ and (unless the algebra is 0,1-dimensional) has infinite volume dimension. It is also acyclic (here $H^0=k1$, the rest are zero). Thus the universal differential calculus it usually too large and interesting by itself.

\subsection{} We define a metric as an element $\cg\in\Omega^1{\bar{\tens}}\Omega^1$ with `inverse' in the sense of a map $\(\ ,\ \):\Omega^1{\bar{\tens}}\Omega^1\to A$ such that
\begin{equation}\label{metric}(\id\tens \(\ ,\omega\))\cg=\omega,\quad (\(\omega,\ )\tens\id)\cg=\omega,\quad\forall \omega\in\Omega^1.\end{equation}
We normally require these maps to be bimodule maps (viewing $\cg:A\to \Omega^1\bar\tens \Omega^1$ by extending its value on 1). This is then equivalent to saying that  the object $\Omega^1$ is left and right self-dual in the monoidal category of $A-A$-bimodules. 

We define a linear (left) connection is a linear map with the property
\begin{equation}\label{conn} \nabla:\Omega^1\to \Omega^1{\bar{\tens}}\Omega^1,\quad \nabla(f\omega)=\extd f{\bar{\tens}}\omega+f\nabla\omega\end{equation}
for all $\omega\in \Omega^1$, $f\in A$. 

When $\Omega^2$ is defined, we define torsion of any connection as 
\begin{equation}\label{tor} T_\nabla:\Omega^1\to\Omega^2,\quad T_\nabla=\wedge\nabla-\extd\end{equation}
where $\wedge:\Omega^1{\bar{\tens}}\Omega^1\to\Omega^2$ is the exterior product. 
Also in this case, we define the curvature of any linear connection as
\begin{equation}\label{curv} R_\nabla:\Omega^1\to \Omega^2{\bar{\tens}}\Omega^1,\quad R_\nabla=(\extd\tens\id-(\wedge\tens\id)(\id\tens\nabla))\nabla.\end{equation}
A similar notion of connection and curvature makes sense on other (usually projective) modules as sections of `vector bundles' but we focus on the cotangent bundle relevant to Riemannian geometry. When there is a metric  say that a connection is `skew metric compatible' or `cotorsion free' if 
\begin{equation}\label{cotor} (\extd\tens\id - (\wedge\tens\id)(\id\tens\nabla))\cg=0.\end{equation}
A torsion free and cotorsion free connection is called `generalised Levi-Civita'. This notion is due to the author. It need not exist but in some examples it does and is unique.

A connection $\nabla$ is called a (left) bimodule connection if there exists a `generalised braiding' $\sigma$ defined in our case by
\begin{equation}\label{braiding} \sigma:\Omega^1{\bar{\tens}}\Omega^1\to \Omega^1{\bar{\tens}}\Omega^1,\quad \nabla(\omega f)=(\nabla\omega)f+\sigma(\omega{\bar{\tens}}\extd f)\end{equation}
for all $\omega\in \Omega^1$, $f\in A$. If $\sigma$ exists then it is a bimodule map. In  this case there is a natural notion of torsion compatibility (which includes the case of torsion free), namely
\begin{equation}\label{omega2} {\rm image}(\id+\sigma)\subseteq\ker\wedge\end{equation}
A bimodule connection naturally extends to  1-1-forms with
\begin{equation}\label{nablag} \nabla\cg =(\nabla\tens\id)\cg+ (\sigma\tens\id)(\id\tens\nabla)\cg,\quad\forall \cg\in \Omega^1\bar\tens\Omega^1\end{equation}
and we say that it is metric compatible on the case of a metric $\cg$  if $\nabla\cg=0$. 
We call a bimodule connection `Levi-Civita' if it is metric compatible and torsion free. This need not exist and we may require a weaker notion of metric and torsion compatible. Bimodule connections were first introduced in \cite{DVM1,DVM2,Mou}.
We refer to \cite{BegMa:star} for an overview of the literature as well as a *-algebra version over $\C$. 

\subsection{} We define a 2nd order differential operator as a linear map $\Delta:A\to A$ such that
\[ \Delta(fg)=(\Delta f)g+f\Delta g+2\<\extd f\bar\tens\extd g\>\]
for some bimodule map $\Omega^1\tens_A\Omega^1\to A$, which will then be uniquely determined. 

Let $(\nabla,\sigma)$ be a bimodule connection and $(\ ,\ )$ a given bimodule map. We do not necessarily require
the connection to be metric compatible. Then
\begin{equation}\label{laplace}\Delta:A\to A,\quad  \Delta= \(\ ,\ \)\nabla \extd. \end{equation}
gives a second order operator which we call the `Laplace-Beltrami' operator associated to $(\nabla,(\ ,\ ))$. Here
\begin{equation}\label{bivec} \<\ ,\ \>= {1\over 2}(\id+\sigma(\ ,\ )).\end{equation}
This follows from the properties of $\nabla,\sigma$ in Section~2.2 and the bimodule property of $(\ ,\ )$. 

Since we know how $\nabla$ acts on 1-1-forms we  likewise extend the Laplace-Beltrami to forms as
\[ \Delta \omega=((\ ,\ )\tens\id)\nabla(\nabla\omega)=((\ ,\ )\tens\id)[(\nabla\tens\id)\nabla\omega+(\sigma\tens\id)(\id\tens\nabla)\nabla\omega]\]
for all $\omega\in \Omega^1$. This time a similar computation gives
\begin{equation}\label{lapder1}
\Delta(f\omega)=(\Delta f)\omega+f\Delta\omega+2\<\ ,\ \>(\extd f\bar\tens\nabla\omega),\quad \forall f\in A,\ \omega\in \Omega^1\end{equation}
for all $\omega,\eta\in\Omega^1$. These construction are new but see \cite{Ma:bh} for an application of such operators.

\subsection{} A differential algebra $(A,\Omega^1,\extd)$ is {\em inner} if there exists $\theta\in\Omega^1$ such that  $[\theta,f]=\extd f$ for all $f\in A$. We use the numbering notation $\sigma_{12}$ for $\sigma$ acting on the first two $\bar\tens$ powers of $\Omega^1$ (and so forth).

 \begin{theorem} Let $A$ be an algebra and $\Omega^1$ an inner differential structure on it.
 \begin{enumerate}\item Bimodule connections $(\nabla,\sigma)$ are in 1-1 correspondence with pairs $(\sigma,\alpha)$
 \[ \sigma:\Omega^1\bar\tens\Omega^1\to \Omega^1\bar\tens\Omega^1,\quad \alpha:\Omega^1\to \Omega^1\bar\tens\Omega^1\]
of bimodule maps and take the form 
 \[ \nabla\omega=\theta\bar\tens\omega-\sigma(\omega\bar\tens\theta)+\alpha\omega.\]
  \item If $\Omega^1$ extends to $\Omega^2$ with $\theta^2=0$ and $\extd\omega=\theta\wedge \omega+\omega\wedge \theta$ for all $\omega\in\Omega^1$ then
 \[ T_\nabla\omega=-\wedge(\id+\sigma)(\omega\bar\tens\theta)+\wedge \alpha\omega\]
 \[ R_\nabla\omega=(\wedge\bar\tens\id)\tilde R_\nabla\omega\]
\[ \tilde R_\nabla\omega=-\sigma_{23}\sigma_{12}(\omega\bar\tens\theta\bar\tens\theta)+(\sigma_{23}(\alpha\tens\id)+(\id\tens\alpha)\sigma)(\omega\bar\tens\theta)-(\id\tens\alpha)\alpha\omega\]
 \item In this case $\nabla$ is torsion free iff it is torsion compatible and $\alpha=0$.
 \item If $\nabla$  is torsion free,  $\sigma(\theta\bar\tens\theta)=\theta\bar\tens\theta$ and $\sigma$ obeys the braid relations, then $R_\nabla=0$.
\end{enumerate}
\end{theorem} 
\proof (1) Let $\sigma, \alpha$ be bimodule maps and $\nabla$ defined from them as stated. Then
\[ \nabla(f\omega)=\theta\bar\tens f\theta-\sigma(f\omega\bar\tens\theta)+\alpha(f\omega)=\theta f\bar\tens \theta-f \sigma(\omega\bar\tens\theta)+f\alpha(\omega)=\extd f\bar\tens\omega+f\nabla\omega\]
\[ \nabla(\omega f)=\theta\bar\tens\omega f-\sigma(\omega\bar\tens f\theta)+\alpha(\omega)f=(\nabla\omega)f+\sigma(\omega\bar\tens\extd f)\]
for all $f\in A$ and $\omega\in \Omega^1$. Hence we have a bimodule connection. Conversely, let $(\nabla,\sigma)$ be a bimodule connection. Then $\sigma$ is a bimodule map and $\nabla^0\omega=\theta\bar\tens\omega-\sigma(\omega\bar\tens\theta)$ is a connection by the above (with $\alpha=0$). Hence $\nabla-\nabla^0$ is a bimodule map, which we take as $\alpha$. (2) We compute
\[ T_\nabla\omega=\wedge\nabla\omega-\extd\omega=\theta\wedge\omega-\wedge\sigma(\omega\bar\tens\theta)+\wedge\alpha\omega-\theta\wedge\omega-\omega\wedge\theta\]
which gives the result as stated. For the curvature we first do the simpler $\alpha=0$ case
\begin{eqnarray*} 
R_{\nabla^0}\omega&=&(\extd\tens\id-(\wedge\tens\id)(\id\tens\nabla^0)(\theta\bar\tens\omega-\sigma(\omega\bar\tens\theta))
\\
&=&\extd\theta\bar\tens\omega-(\extd\tens\id)\sigma(\omega\bar\tens\theta)-\theta^2\bar\tens\omega+\theta\wedge\sigma(\omega\bar\tens\theta)+(\wedge\tens\id)(\id\tens\nabla^0)\sigma(\omega\bar\tens\theta)\\
&=&-\sigma_1\wedge\theta\bar\tens\sigma_2+(\wedge\tens\id)(\id\tens\nabla^0)\sigma(\omega\bar\tens\theta)=-(\wedge\tens\id)\sigma_{23}\sigma_{12}(\omega\bar\tens\theta\bar\tens\theta)
\end{eqnarray*}
on using the definition of $\nabla^0$. Here $\sigma_1\bar\tens\sigma_2:=\sigma(\omega\bar\tens\theta)$ is a notation. Then
\begin{eqnarray*}
R_\nabla\omega&=&(\extd\tens\id-(\wedge\tens\id)(\id\tens\nabla))(\theta\bar\tens\omega-\sigma(\omega\bar\tens\theta)+\alpha\omega)\\
&=&R_{\nabla^0}\omega+(\extd\tens\id)\alpha\omega-\alpha_1\wedge(\theta\bar\tens\alpha_2-\sigma(\alpha_2\bar\tens\theta)+\alpha\alpha_2)-\theta\wedge\alpha\omega+\sigma_1\wedge\alpha\sigma_2
\end{eqnarray*}
where $\alpha\omega:=\alpha_1\bar\tens\alpha_2$ is a notation. Using the form of $\extd$ and cancelling two terms we arrive at the expression stated. (3) If (\ref{omega2}) holds and $\alpha=0$ then we see from part (2) that $T_\nabla=0$. Conversely, torsion free is a special case of torsion-compatibility. (4) By part (3) the assumption entails torsion compatible and  $\alpha=0$. If $\sigma$ preserves $\theta\bar\tens\theta$ (a sufficient condition for our assumption that $\theta^2=0$) and the braid relations (or `Yang-Baxter equations') $\sigma_{23}\sigma_{12}\sigma_{23}=\sigma_{12}\sigma_{23}\sigma_{12}$ hold,  then
\[ R_\nabla\omega=-(\wedge\tens\id)\sigma_{23}\sigma_{12}(\omega\bar\tens\theta\bar\tens\theta)=-(\wedge\tens
\id)\sigma_{23}\sigma_{12}\sigma_{23}(\omega\bar\tens\theta\bar\tens\theta)\]
\[=-(\wedge\tens\id)\sigma_{12}\sigma_{23}\sigma_{12}(\omega\bar\tens\theta\bar\tens\theta)[=(\wedge\tens\id)\sigma_{23}\sigma_{12}(\omega\bar\tens\theta\bar\tens\theta)=-R_\nabla\omega\]
where we used (\ref{omega2}). If the characteristic  is not 2 we conclude that $R_\nabla=0$. \endproof
 
We will be mainly interested in this simplest case $\alpha=0$, which we call the  bimodule connection associated to a bimodule map $\sigma$. 

\begin{lemma}  Let $A$ be an algebra, $\Omega^1$ an inner differential structure and $\sigma$ a bimodule map on $\Omega^1\bar\tens\Omega^1$.
\begin{enumerate}\item If $\sigma(\theta\bar\tens\theta)=\theta\bar\tens\theta$ then we meet the conditions on $\Omega^2$ in the theorem by defining
  \[ \Omega^2=\Omega^1\bar\tens\Omega^1/\ker(\id-\sigma) \]
  \item In this case the bimodule connection corresponding to $(\sigma,0)$ has $R_\nabla\omega=0$ iff $\sigma$ obeys the braid relations on $\omega\bar\tens\theta\bar\tens\theta$. \end{enumerate}
\end{lemma}
\proof (1) Under the assumption $\theta\tens\theta\in \ker (\id-\sigma)$ we have $\theta^2=0$. We define $\extd\omega=\theta\wedge\omega+\omega\wedge\theta$ which is then necessarily a graded derivation with respect to products by functions. We have to check that $\extd\extd f=0$ but this follows provided $\theta^2$ commutes with functions, in particular if it vanishes. (2) Let $\tilde R_\nabla$ be the expression in Theorem~2.1 so that $R_\nabla\omega=0$ iff $\tilde R_\nabla\omega\in \ker\wedge$ i.e. iff $\tilde R_\nabla\in\ker(\id-\sigma)$. This happens iff
\[ \sigma_{12}\sigma_{23}\sigma_{12}(\omega\bar\tens\theta\bar\tens\theta)=\sigma_{23}\sigma_{12}(\omega\bar\tens\theta\bar\tens\theta)\]
but the right hand side is also $\sigma_{23}\sigma_{12}\sigma_{23}(\omega\bar\tens\theta\bar\tens\theta)$. \endproof

We call this construction the `canonical prolongation' to degree 2 associated to a bimodule map $\sigma$ that leaves $\theta\bar\tens\theta$ invariant. At least when $\sigma$ obeys the braid or Yang-Baxter relations there is a similar prolongation to all degrees as the tensor algebra on $\Omega^1$ modulo the kernel of `braided antisymmetrizers' $[n,-\sigma]!$ in degree $n$. We use here a canonical lift of $S_n$ to the braid group $B_n$ provided by  braided factorials \cite{Ma:book} in which simple reflections in a reduced expression are replaced by $\sigma$ in the corresponding position (cf. a construction in a quantum group context due to Woronowicz). The exterior derivative is $\extd=[\theta,\ \}$ (a graded commutator) and we obtain a differential graded algebra (the noncommutative de Rham complex) and its associated cohomology. 

\subsection{} We define a `non-standard-Ricci curvature'  as follows. Classically, one lifts the 2-form value of $R_\nabla$ to values in $\Omega^1\bar\tens\Omega^1$ and then takes a trace. This lift requires more structure but in our case above we have an obvious but non-classical candidate $\tilde R_\nabla$ featuring in Theorem~2.1 and such that $R_\nabla=(\wedge\tens\id)\tilde R_\nabla$. In place of the trace it is convenient to use the metric and inverse metric  leading to
\begin{equation}\label{ricci} S_{\nabla,\cg}=\(\ ,\ \)_{12}(\id\tens\tilde R_\nabla)\cg\in \Omega^1\bar\tens\Omega^1.\end{equation}
When $\Omega^1$ is finitely generated and projective one can replace the metric and inverse metric here by a trace, as well as consider more classical lifts of $R_\nabla$. There is not yet a definitive notion of Ricci tensor at the level of differential algebras.

\section{Metrics and connections on graphs}          

Here we translate some of the algebraic Riemannian geometry in Section 2 to the   level of graphs, starting with a clear `dictionary'. We then construct metrics and connections $\nabla$ including a canonical metric and connection on any graph  and a metric compatible  `Maurer-Cartan' connection on any Cayley graph. 

\subsection{Functorial correspondence}

Let $A=k(V)$ be the algebra of functions on a finite set $V$ with pointwise product. From the Pierce decompisiton of $A$ (i.e. by considering the Kronecker delta function projectors $\{\delta_x\ |\ x\in V\}$) it is easy to see that the possible $\Omega^1$ are in 1-1 correspondence with subsets $E\subseteq V\times V\setminus{\rm diagonal}$, i.e. to digraph structures $(V,E)$ where $E$ denotes the set of arrows. We write $x\to y$ iff $(x,y)\in E$ according to this identification. By a  digraph throughout the paper we mean that 
 are no arrows from a vertex to itself and at most one arrow in each direction between distinct vertices. This observation is not new, see \cite{Connes}. Explicitly, 
 \[ \Omega^1=k E,\quad f\omega_{x \to y}=f(x)\omega_{x\to y},\quad \omega_{x\to y}f=\omega_{x\to y}f(y),\quad \extd f=\sum_{x\to y\in E}(f(y)-f(x))\omega_{x\to y}\]
for all functions $f\in k(V)$ and basis $\{\omega_{x \to y}\}$ labelled by the arrows $x\to y\in E$.

Note that 
\[ \omega_{x\to y}\bar\tens\omega_{z\to w}=\delta_{y,z}\omega_{x\to y}\bar\tens\omega_{y\to w}\]
where $\bar\tens$ is the tensor product over $k(V)$. Similarly for higher $\bar\tens$ products.  We also note that
\[ \extd\delta_x=\sum_{y: y\to x}\omega_{y\to x}-\sum_{y:x\to y}\omega_{x\to y},\quad \delta_x\extd\delta_y=\begin{cases} -\sum_{z:x\to z}\omega_{x
\to z}& x=y\\ \omega_{x\to y}& x\to y\\ 0 &{\rm else}\end{cases}. \]

Possibly less well-known but important in what follows: the differential algebra $\Omega^1(V,E)$ associated to a digraph $(V,E)$ is always inner, with
\[ \theta=\sum_{x\to y\in E}\omega_{x\to y}\]
This is immediate from the bimodule relations stated.

Let {\rm Digraph} be the category of finite directed graphs with the following notion of morphisms. In fact there are different notions in the literature\cite{HelNes} and the most usual of sending vertices to vertices and arrows to arrows is not quite what we want, partly because we have decided to work with simple digraphs in which there are no self-arrows. In our case we define a morphism of digraphs $(W,F)\to (V,E)$ as a map 
\[ \psi:W\to V,\quad {s.t.}\quad   \forall w\to z\in F, \quad  \psi(w)=\psi(z)\quad{\rm or}\quad \psi(w)\to\psi(z)\in E.\]
 A cleaner point of view is to equivalently consider our digraphs as `extended-simple' digraphs where every vertex has a self-arrow understood. The corresponding map between such extended-simple digraphs is simply a map of vertices that sends arrows to arrows. In fact an extended-simple digraph is a category with objects the vertices and morphisms the arrows, and a map between extended-simple digraphs is a functor. From this point of view Digraph actually forms a 2-category with morphisms between $\psi,\psi'$ the natural transformations. Explicitly $\psi\Rightarrow\psi'$ if for all $w\in W$ either  $\psi(w)=\psi(w')$ or $\psi(w)\to \psi'(w)\in E$ such that 
\[ \begin{array}{ccc}\psi(w)&\to &\psi'(w)\\ \downarrow& & \downarrow\\ \psi(z)&\to &\psi'(z)\end{array}  \]
for all $w\to z\in F$ in the generic case and similarly triangles if $\psi(w)=\psi'(w)$ or $\psi(z)=\psi'(z)$. Again we need not single out the equality cases if we work with extended-simple digraphs. In other words, for every two fixed digraphs the set of morphisms from one to the other is itself a digraph in a natural way.  Let ${\rm DiffAlg}$ be the category of differential algebras differentiable algebra maps in the sense of (\ref{moralg}) as morphisms. 

\begin{theorem} The association above extends to a full and faithful functor ${\rm Digraph}\to {\rm DiffAlg}$. \end{theorem}
\proof We first check that the map $(V,E)\mapsto (k(V),\Omega^1(V,E),\extd)$ extends to a functor. A morphism $(W,F)\to (V,E)$ of digraphs has been described above as a type of map $\psi:W\to V$. This induces a map $\phi:k(E)\to k(F)$ by pull-back and a further map 
\begin{equation}\label{phi*} \phi_*:\Omega^1(V,E)\to \Omega^1(W,F),\quad \phi_*(\omega_{x\to y})=\sum_{ w\to z\atop  w\in \psi^{-1}(x),\ z\in\psi^{-1}(y)}\omega_{w\to z}\end{equation}
We check that this is a bimodule map: if $f\in k(E)$, we have $\phi_*(f \omega_{x\to y})= f(x)\phi_*(\omega_{x\to y})$, computing on the other side of the equality $\phi(f)\phi_*(\omega_{x\to y})$ has a factor in each term $\phi(f)(w)=f(\psi(w))=f(x)$ as required. Similarly for the action from the right. Next we show that $\phi_*(\extd f)=\extd\phi(f)$ but since both calculi are inner it suffices to show that $\phi_*(\theta)=\theta$ for the two calculi. Indeed
\[ \phi_*(\sum_{x\to y}\omega_{x\to y})= \sum_{x\to y}\sum_{ w\to z\atop  w\in \psi^{-1}(x),\ z\in\psi^{-1}(y)}\omega_{w\to z}\]
but any $w\to z$ appears on the sum on the right since $\psi(w)\to\psi(z)$ is an arrow. Also no two arrows in $F$ can appear twice in the right hand sum because for each $x\to y$ the arrows are between fibers of $\psi$ and hence cannot overlap with other summands as $x,y$ change. Hence the right hand side a sum over all arrows of $F$, is the inner element $\theta$ for the differential structure associated to $(W,F)$. This gives our functor on morphisms. 

Suppose $\psi,\psi':(W,F)\to (V,E)$ are distinct digraph morphisms. If they are the same on the vertex sets then they are the same since the action on arrows is determined, hence they differ on the vertex sets and the induced $\phi,\phi':k(V)\to k(W)$ differ as algebra maps between the given algebras. So the functor is clearly faithful. Now suppose $(\phi,\phi_*)$ is a morphism of differential algebras with $\phi_*:\Omega^1(V,E)\to \Omega^1(W,F)$ and $\phi:k(V)\to k(W)$. The latter being an algebra map must have the form $\phi(f)(w)=f(\psi(w))$ for some set map $\psi:W\to V$. This is obvious from the Pierce decomposition but to see it directly let $\phi(\delta_x)=\sum_{w\in W}\phi^x_w\delta_w$ and note that since $\phi(1)=\phi(\sum_x\delta_x)=1=\sum_w\delta_w$ we have $\sum_{x\in V}\phi^x_w=1$ for all $w$. Meanwhile $\phi(\delta_x)=\phi(\delta^2_x)=\sum_{w,z}\phi^x_w\phi^x_z\delta_{w,z}\delta_w=\sum_{w}\phi^x_w\phi^x_w\delta_w$ from which we conclude that $\phi^x_w\in\{0,1\}$. Hence for each $w$ there is precisely one $x$ where $\phi^x_w\ne 0$ and we call this $\psi(w)$. Since we know that $\phi$ is differentiable the map $\phi_*$ is uniquely determined and must therefore be the map induced by $\psi$. Alternatively we have seen that  $\phi(\delta_x)=\sum_{w\in\psi^{-1}(x)}\delta_w$ and hence $\phi_*(\delta_x\extd\delta_y)=\phi(\delta_x)\extd\phi(\delta_y)$ comes out as the stated form of $\phi_*(\omega_{x\to y})$ when $x\to y$.  When $x\ne y$ and not $ x\to y$, we must get zero from this expression i.e. there do not exist $w\in\psi^{-1}(x)$ and $z\in \psi^{-1}(y)$ such that $w\to z$. The contrapositive of this is that if $w\to z$ then $\psi(w)\to\psi(z)$ or $\psi(w)=\psi(z)$.

\endproof

In particular, two digraphs are isomorphic iff their associated differential algebras which means that we can transfer ideas from one category to the other systemmatically. We also know a little more -- there are no objects of the form $(k(V),\Omega^1,\extd)$ in ${\rm DiffAlg}$ other than in the image of this functor. We illustrate some elements of this correspondence. For example, the universal differential structure on $k(V)$ corresponds to the complete digraph on the vertex set $V$ and its corresponding property is that every digraph on $V$ is a subdigraph. 

\begin{proposition} Let be a digraph $(V,E)$. The corresponding differential algebra $\Omega^1(V,E)$ is:
\begin{enumerate}\item  connected {\em iff} the digraph is weakly connected (i.e. the underlying graph is connected)
\item  left(right) parallelizable of dimension $n$ {\em iff} the digraph is $n$-regular for outgoing (incoming) arrows.
\item  both left and right parallelizable {\em iff} the digraph has an equal and constant number of in and out arrows at every vertex. In this case the left and right dimensions are the same and there is a simultaneous left and right basis.
\end{enumerate}
\end{proposition}
\proof (1) From the form of $\extd$, if $\extd f=0$ then $f(x)-f(y)=0$ for all $x\to y$. Hence if the graph is weakly connected we conclude that $f$ is constant. Conversely, if the differential structure is connected. We recall that a digraph is bidirected if every arrow has a reverse arrow going the other way. 
(2)  If $\{\omega_i\}$ are  a global left basis then $\{\delta_x\omega_i\}$ are necessarily a basis of $\delta_x\Omega^1$ for each $x$. Here, if $\delta_x\omega\in\delta_x\Omega^1$ we can write $\omega=\sum f_i\omega_i$ for some $f_i$ hence $\delta_x\omega= \sum f_i(x)\delta_x\omega_i$ so they span. And if $\sum\lambda_i\delta_x\omega_i=0$ we let $f_i(y)=\lambda_i\delta_{x,y}$ for all $y\in X$ and conclude that $\delta_y\sum_if_i\omega_i=0$ for all $y$, hence $\sum_i f_i\omega_i=0$. This then requires all the $\lambda_i=0$. Hence $\dim(\delta_x\Omega^1)=n$ but a basis of $\delta_x\Omega^1$ is $\{\omega_{x\to y}\ |\ y\}$ hence there are precisely $n$ arrow out of each vertex.  Conversely, suppose the graph is outgoing $n$-regular and choose a colouring of the arrows out of each vertex by  $i\in\{1,\cdots, n\}$ (so that the arrows are enumerated). Let 
\[ \omega_i=\sum_{x{\buildrel i \over \to }y}\omega_{x\to y}\]
then one can verify that this is a basis. Here every $\omega\in \Omega^1$ is a $k$ linear combination of the $\omega_{x\to y}$ and so we can write $\omega=\sum_i f_i\omega_i$ where $f_i$ are defined so that $f_i(x)$ is the coefficient of $\omega_{x{\buildrel i\over \to }y}$. (For  a fixed $i$ each $x\in X$ occurs just once as source of an arrow and there is a unique $y$ as its head). Clearly if $\sum_i f_i\omega_i=0$ then the $f_i=0$ as the $\{\omega_{x\to y}\}$ are a basis over $k$. The right handed result is analogous. (3) Finally, as the sum of the in-degrees and the sum of the out-degrees are equal (being the number of arrows), if the differential structure is both left and right parallelizable then the left and right dimensions must be the same and hence the in and out valencies the same (say $n$) at every vertex. However, we can say more. By Hall's marriage theorem\cite{Hall} every digraph of this type is a union of $n$ 1-difactors. To see this, double the vertices to give a bipartite graph with the original arrows as edges from one `outgoing' set of vertices to the other `incoming' set of vertices. The marriage theorem gives a  bijection of one to the other via a subset of the edges and this traces out a 1-difactor. Deleting these edges we have a similar situation for $n-1$ vertex in and out degrees, i.e. the proof is by indiction.  We can then assign a colour to each of these 1-difactors and hence colour the original digraph so that at each vertex every colour occurs precisely once as an out arrow and precisely once as an in arrow. In this case the $\{\omega_i\}$ above will have the desired basis property as in part (2) but from both the left and the right at the same time. \endproof

\begin{example} The graph
\[ \includegraphics[scale=.5]{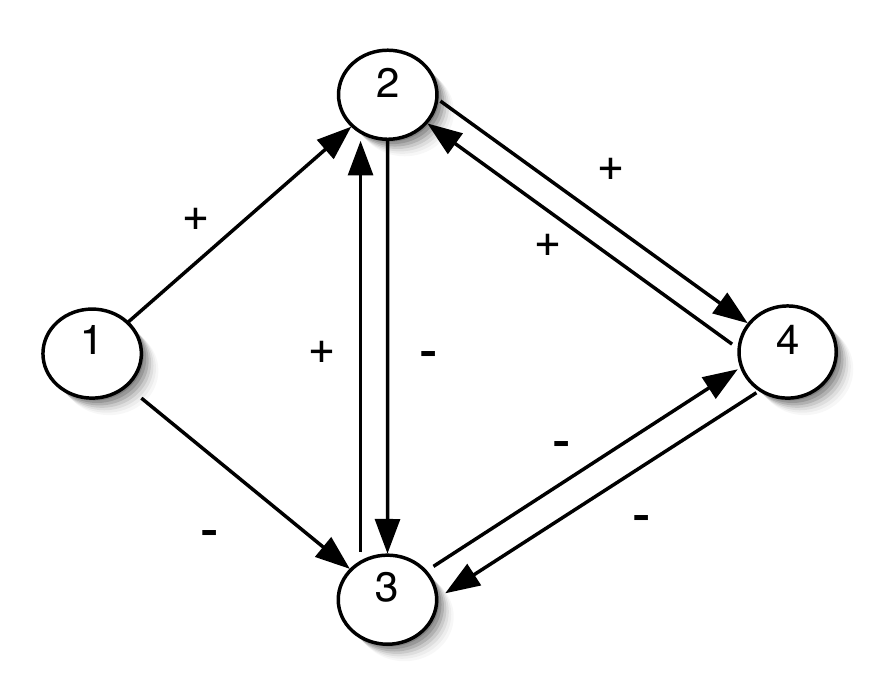}\]
is left regular with left-dimension 2 but not right-regular at all.  The edge labels shown are an example of an edge colouring by $\{+,-\}$ with the associated left global basis. 
\[ \omega^+=\omega_{1\to 2}+\omega_{2\to 4}+\omega_{4\to 2}+\omega_{3\to 2},\quad\omega^-=\omega_{1\to 3}+\omega_{2\to 3}+\omega_{3\to 4}+\omega_{4\to 3}\]
Observe that each vertex occurs precisely once as a tail of an arrow in $\omega^+$ and similarly in $\omega^-$, leading to the left basis property.
\end{example}

The maximal prolongation to $\Omega^2$ depends only on $\Omega^1$ which combined with the functoriality in Theorem~3.1 means that the noncommutative de Rham cohomology $H^1(k(V),\Omega(V,E),\extd)$ is an invariant of digraphs. The same will apply to `natural' choices of $\Omega^2$ depending on the data used. An existing description of the maximal prolongation of $\Omega^1(V,E)$ and hence of the data for $\Omega^2(V,E)$ is in \cite{BrzMa}.  We shall take a fresh approach where $\Omega^2$ is defined by a connection using Lemma~2.2. 

\subsection{General metrics and connections on bidirected graphs}

We recall that a digraph is bidirected if every arrow has a reverse arrow going the other way. 

\begin{proposition} The differential structure associated to a digraph $(V,E)$ admits a metric {\em iff} it is bidirected. The metric then takes the form
\[ \cg=\sum_{x\to y\in E}{1\over g_{x\to y}}\omega_{x\to y}\bar\tens\omega_{y\to x},\quad \(\omega_{x\to y},\omega_{y'\to x'}\)=g_{y\to x}\delta_{x,x'}\delta_{y,y'}\delta_x\]
for arbitrary $\{g_{x\to y}\in k^*\}$ associated to each arrow.
\end{proposition}
\proof $\Omega^1\bar\tens\Omega^1$ is spanned by elements of the form $\omega_{x\to y}\bar\tens\omega_{y\to z}$ in order for the tensor product to be well-defined over $k(V)$ (one can insert $\delta_y$ in the middle and would otherwise get zero). In order to be a bimodule map we need $\cg$ to commute with functions and this forces $z=x$. This forces us to sum only over arrows which have reverse arrows, giving the form of $\cg$ stated but with the sum over such arrows. Similarly for the form of $(\ ,\ )$ which will be zero of the relevant arrows are not reversiible. If $x'\to y'$ is a reversible arrow, we verify (\ref{metric}),
\begin{eqnarray*} \sum_{x\to y}{1 \over g_{x\to y}}\(\omega_{x'\to y'},\omega_{x\to y}\)\omega_{y\to x}=\sum_{x\to y}{g_{y'\to x'}\over g_{x\to y}}\delta_{y',x}\delta_{x',y}\delta_{x'}\omega_{y\to x}=\omega_{x'\to y'}\\
\sum_{x\to y}{1 \over g_{x\to y}}\omega_{x\to y}\(\omega_{y\to x},\omega_{x'\to y'}\)=\sum_{x\to y}{g_{x\to y}\over g_{x\to y}}\omega_{x\to y}\delta_{x,x'}\delta_{y,y'}\delta_y=\omega_{x'\to y'}
\end{eqnarray*} using the left and right bimodule structures to evaluate the delta-functions. However, we will not be able to obey (\ref{metric}) for $\omega_{x'\to y'}$ associated to arrows that are not reversisble, so we will need the digraph to be bidirectional for it to hold on all of $\Omega^1$. \endproof

If we are given an undirected graph we can view it as a bidirected digraph by assigning arrows in both directions to each edge, but note that there are still two basis vectors $\omega_{x\to y}, \omega_{y\to x}\in\Omega^1$  associated to the edge. In view of Proposition~3.4,  throughout the rest of the section, in fact the rest of the paper, we will focus on undirected graphs which we view as bidirected. We will still indicate arrows in quantifiers and conditions in order to be definite and because the direction of arrows is typically connected with the computation being made, but it should be remembered that if $x\to y$ is an arrow, so is $y\to x$ so that the direction in a quantifier such as $\forall x\to y$ is imaterial. 
We similarly analyse the possible connections in this bidirected case.

\begin{lemma} Let $\Gamma$ be an undirected graph and $\Omega^1$ its bidirected differential structure. Bimodule connections on $\Omega^1$ are defined by numerical data $\sigma,\alpha$ where
\[ \sigma^{x,y,z}_{w}\in k,\quad\forall\quad  \begin{array}{rcl}x\to & y&\to z\\ \searrow & w & \nearrow\end{array},\quad\quad  \alpha^{x,y}_w\in k,\quad \forall \quad \begin{array}{rcl}x&\longrightarrow&y\\ &\searrow w\nearrow& \end{array}\]
and
\begin{eqnarray*} \nabla\omega_{x\to y}&=&\sum_{z:z\to x\to y}\omega_{z\to x}\bar\tens\omega_{x\to y}\\
&&-\sum_{w,z:\ \begin{array}{rcl}x\to & y&\to z\\ \searrow & w & \nearrow\end{array}}\sigma^{x,y,z}_w\omega_{x\to w}\bar\tens\omega_{w\to z}+\sum_{w:\ \begin{array}{rcl}x&\longrightarrow&y\\ &\searrow w\nearrow& \end{array}}\alpha^{x,y}_w\omega_{x\to w}\bar\tens\omega_{w\to y}\end{eqnarray*}
\end{lemma}
\proof The sums here are over all $w,z\in V$ or all $w\in V$ for which there exist edges in the graph as shown. The directions of the arrows are irrelevant to these criteria since all our edges are bidirectional but included for clarity. We know from Theorem~2.1 that a bimodule connection has the form $\nabla\omega=\theta\bar\tens\omega-\sigma(\omega\bar\tens\theta)+\alpha\omega$ for all $\omega\in\Omega^1$. The only thing we need to note is that  a bimodule map $\sigma:\Omega^1\bar\tens\Omega^1\to \Omega^1\bar\tens\Omega^1$ necessarily has the form
\[\sigma(\omega_{x\to y}\bar\tens\omega_{y\to z})=\sum_{w:\ \begin{array}{rcl}x\to & y&\to z\\ \searrow & w & \nearrow\end{array}}\sigma^{x,y,z}_w\omega_{x\to w}\bar\tens\omega_{w\to z}.\]
for some coefficients $\sigma^{x,y,z}_w$. Similarly a bimodule map $\alpha:\Omega^1\to \Omega^1\bar\tens\Omega^1$ necessarily has the form
\[ \alpha \omega_{x\to y}=\sum_{w:\ \begin{array}{rcl}x&\longrightarrow&y\\ &\searrow w\nearrow& \end{array}}\alpha^{x,y}_w\omega_{x\to w}\bar\tens\omega_{w\to y}\]
for some coefficients $\alpha^{x,y}_w$. 
\endproof
	
 In passing, we write out the braid or Yang-Baxter relations for $\sigma$.

\begin{proposition} Using notations as above, a bimodule map  $\sigma:\Omega^1\bar\tens\Omega^1\to \Omega^1\bar\tens\Omega^1$ obeys the braid relations iff
\[\forall \begin{array}{rcl} &\nearrow y\searrow&\\ x& &z\\ \downarrow& &\downarrow\\ t& &w\\ &\searrow p \nearrow\end{array},\quad \sum_{s}\sigma^{y,z,w}_s\sigma^{x,y,s}_t\sigma^{t,s,w}_p=\sum_{s}\sigma^{x,y,z}_s\sigma^{s,z,w}_p\sigma^{x,s,p}_t  \]
where the sum is over $s$ and the edges are according to the respective sides of the diagram
\[ \includegraphics[scale=.35]{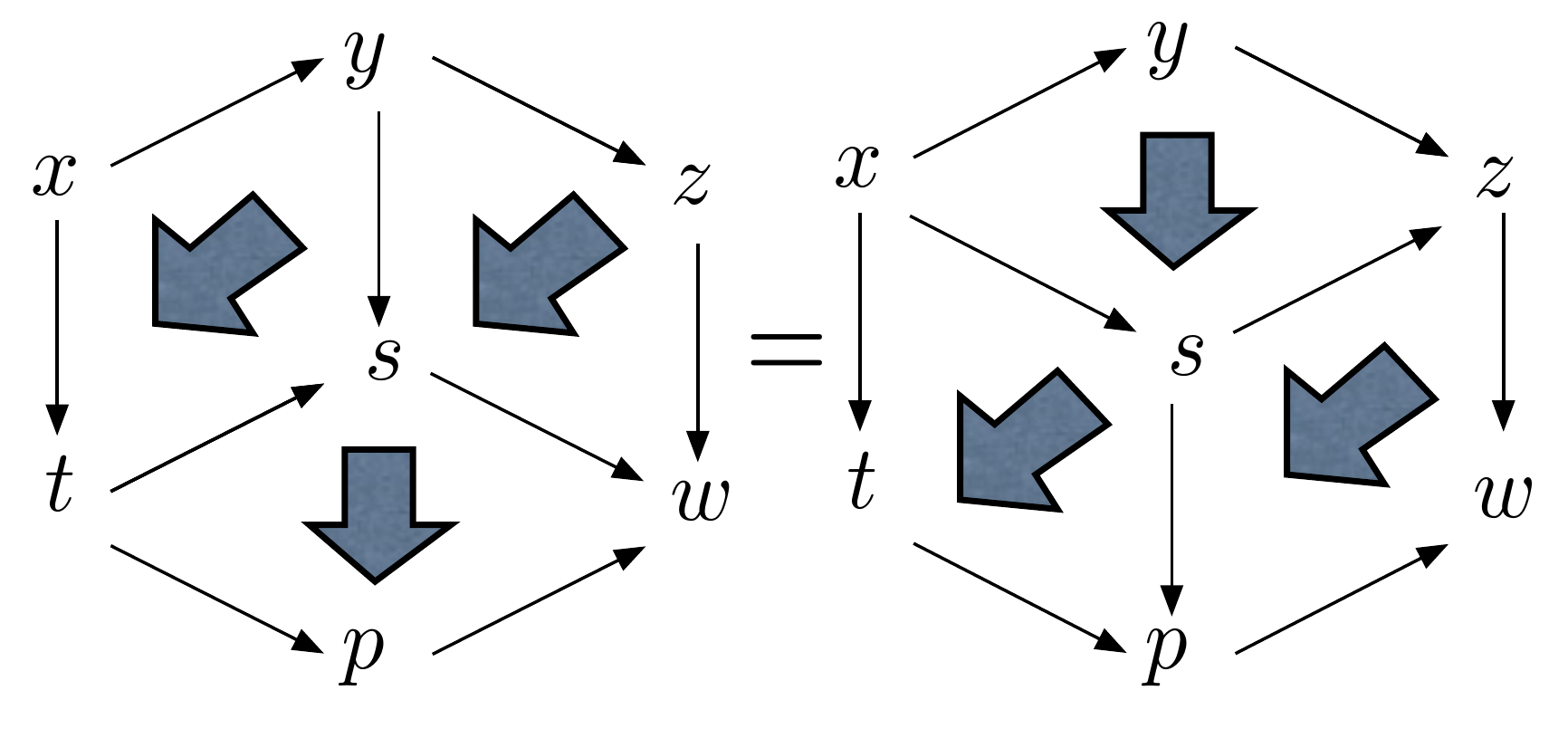}\]\
The diagram also depicts $\sigma$, as a solid arrow. 
\end{proposition}
\proof This is immediate from the definition of $\sigma$. \endproof
One can view the braid relations here as equality of  two ways to `surface transport' from the boundary $x\to y\to z\to w$ to the boundary $x\to t\to p\to w$ by application of $\sigma$, namely  round the front faces or round the back faces of a cube as shown. On either side, all possible values of the internal vertex are summed over. 

\goodbreak
\begin{proposition} {\  }
\begin{enumerate} \item For any connection $(\sigma,\alpha)$ and metric $\cg$ the Laplace-Beltrami operator  (\ref{laplace}) is the weighted graph Laplacian
\[ (\Delta f)(x)=\sum_{y:x\to y}( f(x)-f(y))\gamma_{x,y},\quad \gamma_{x,y}=g_{y\to x}+\sum_{w: x\to w} g_{w\to x}\sigma^{x,y,x}_w\]
\item $(\ ,\ )\sigma=(\ ,\ )$  {\em iff} $\sum_{w: x\to w} g_{w\to x}\sigma^{x,y,x}_w=g_{y\to x}$ and in this case $\gamma_{x,y}=2g_{y\to x}$. \end{enumerate}
\end{proposition}
\proof We compute $\Delta f=\sum_{x\to y}(f(y)-f(x))\(\ ,\ \) \nabla\omega_{x\to y}$ where
\begin{eqnarray*} \(\ ,\ \) \nabla\omega_{x\to y} &=&(\omega_{y\to x},\omega_{x\to y})-\sum_{ \begin{array}{rcl}x\to & y&\to x\\ \searrow & w & \nearrow\end{array}}\sigma^{x,y,z}_w(\omega_{x\to w},\omega_{w\to x}) \\
&=&g_{x\to y}\delta_y -\sum_{w:x\to w}\ g_{w\to x}\delta_x\sigma^{x,y,x}_w\end{eqnarray*}
where the contribution from $\alpha$ vanishes since we do not have $x=y$ which would be needed for a non-zero value of $(\ ,\ )$. We use its form from Proposition~3.4 and the form of $\nabla$ in Lemma~3.5. This gives the expression stated after a change of variables in one of the terms. The second part is clear from the computation in the proof. \endproof

We see that noncommutative geometry reproduces  the usual weighted graph-theory Laplacian with weight determined form the metric and connection. We see that if $(\ ,\ )\sigma=(\ ,\ )$ which is a natural simplifying condition in the general theory of Section~2.3, the weights are the metric coefficients, just as in classical geometry. 

We now consider the question of full metric compatibility but for simplicity we specialise to the case $\alpha=0$. The general case is no harder but the formulae are more complicated and we will not be needing them.

\begin{lemma} A connection defined by $(\sigma,0)$ is metric compatible for the metric $\cg$ iff
\[ \forall \quad  \begin{array}{rcl}x&\longrightarrow&z\\ \downarrow & &\uparrow\\ v&\to&w \end{array},\quad \sum_{y:x\from y\to w}\sigma^{y,x,z}_w\sigma^{x,y,w}_v{g_{z\to w}\over g_{x\to y}}=\delta^z_v\]
We depict this equation in diagrammatic form
\[ \includegraphics[scale=.7]{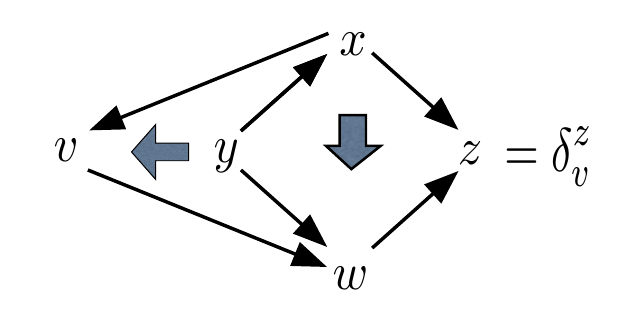}\]
where there is a sum over $y$. 
\end{lemma}
\proof We proceed from the definitions of $\nabla, \cg$ and $\nabla\cg$ in Proposition~3.4 and Lemma~3.5,
\begin{eqnarray*}
\nabla\cg&&= \sum_{z\to x\to y}{1\over g_{x\to y}}\omega_{z\to x}\bar\tens\omega_{x\to y}\bar\tens\omega_{y\to x}-\kern-15pt \sum_{\begin{array}{rcl}y\to & x &\to z\\ \searrow & w & \nearrow\end{array}}\kern-10pt{ \sigma^{y,x,z}_w\over g_{x\to y}}\sigma(\omega_{x\to y}\bar\tens\omega_{y\to w})\bar\tens\omega_{w\to z}\\
&&= \sum_{x\to z\to w}{1\over g_{z\to w}}\omega_{x\to z}\bar\tens\omega_{z\to w}\bar\tens\omega_{w\to z} -\kern-15pt \sum_{\begin{array}{rcl}& x & \\ &\nearrow\  \downarrow\  \searrow&\\ y& z& v\\ &\searrow\ \uparrow\ \swarrow&\\ &w&\end{array}}\kern-10pt {\sigma^{x,y,w}_v\sigma^{y,x,z}_w\over g_{x\to y}}\omega_{x\to v}\bar\tens\omega_{v\to w}\bar\tens\omega_{w\to z}\end{eqnarray*}
where we do not show two terms that immediately vanish due to the $\bar\tens$. In the second equality we have relabelled some of the variables being summed over. The diagrammatic form indicates  apply $\sigma$. Note that $y$ is a summed variable and that its arrows are not significant. \endproof

\subsection{Canonical Euclidean metric and canonical connection}

In this section we consider particularly the `Euclidean metric' 
\[ \cg= \sum_{x\to y}\omega_{x\to y}\bar\tens \omega_{y\to x},\quad \(\omega_{x\to y},\omega_{y'\to x'}\)=\delta_{x,x'}\delta_{y,y'}\delta_x\]
where there are no weights factors $g_{x\to y}$. We also consider what we call the  `canonical connection' 
\[ \nabla\omega_{x\to y}=\sum_{z\to x}\omega_{z\to x}\bar\tens \omega_{x\to y}-\sum_{y\to z}\omega_{x\to y}\bar\tens\omega_{y\to z},\quad  \sigma=\id\]
which corresponds to $\sigma=\id$ and $\alpha=0$ in the general theory of Section~2.4 and to  $\sigma^{x,y,z}_w=\delta^y_w$ and $\alpha^{x,y}_w=0$ in the explicit
notation of Section~3.2. We call these canonical because they do not involve any weights and therefore make sense over any field. As such they would also apply in
a scheme theory setting.

The associated Laplace-Beltrami operator to these two objects is $\Delta= 2L$ where $L$ is the canonical unweighted graph Laplacian
\[ (Lf)(x)=\sum_{d(x,y)=1}f(x)-f(y)=\sum_{y:x\to y} f(x)-f(y)\]
where $d(x,y)$ is the minimum number of edges between $x$ and $y$ (so the sum is over nearest neighbours of $x$. As part of the general theory, we deduce that this obeys 
\[ L(fg)=(Lf)g+f(Lg)+(\extd f,\extd g),\quad \forall f,g\in k(V)\]
Also in the case of the Euclidean metric, it is easy to see that 
\[ (\theta,\theta)=\deg\]
is the degree function, where $\deg(x)$ is the number of undirected edges at the vertex $x$. 

Although both appear canonical in a graph theory context, note that $\nabla$ here is not metric compatible with $\cg$ here (it does not need to be in our framework). However,  the natural space of 2-forms provided by Theorem~2.1 for $\sigma=\id$  is $\Omega^2=0$, and says that the canonical  geometry here is in some sense 1-dimensional. Because of this, the connection is trivially torsion-free and trivially cotorsion free so one could say that it is `generalised Levi-Civita' as explained in Section~2.2, just because $\Omega^2=0$. Also clearly the curvature is zero and we note also that $\sigma=\id$ trivially obeys the braid or Yang-Baxter equations. 
The connection $\nabla$ is, however,  strange compared to classical geometry, as the following example illustrates.

\begin{example} In the case of the 1-dimensional line graph we identify the vertex set $V=\Z$. Then our canonical non-classical connection is clearly
\[ \nabla\omega_{x\to x+1}=\omega_{x-1\to x}\bar\tens\omega_{x\to x+1}+\omega_{x+1\to x}\bar\tens\omega_{x\to x+1}-\omega_{x\to x+1}\bar\tens\omega_{x+1\to x}-\omega_{x\to x+1}\bar\tens\omega_{x+1\to x+2}\]
\[ \nabla\omega_{x+1\to x}=\omega_{x\to x+1}\bar\tens\omega_{x+1\to x}+\omega_{x+2\to x+1}\bar\tens\omega_{x+1\to x}-\omega_{x+1\to x}\bar\tens\omega_{x\to x-1}-\omega_{x+1\to x}\bar\tens\omega_{x\to x+1}.\]
In order to understand this, we note (and we will recall this later) that the group structure of $\Z$ implies that $\Omega^1$ is a free module over $k(V)$ with basis the left-invariant 1-forms
\[ e_+=\sum_x \omega_{x\to x+1},\quad e_-=\sum_x\omega_{x+1\to x}\]
After elementary computations we find
\[ \cg=e_+\bar\tens e_-+e_-\bar\tens e_+,\quad \nabla e_+=-\nabla e_-=e_-\bar\tens e_+-e_+\bar\tens e_-.\]
We can understand this more `geometrically' by considering the linear coordinate function $X(x)=x$.  Then $\extd X=e_+- e_-$ and we also have $\theta=e_++e_-$. In these terms 
\[ \cg={1\over 2}(\theta\bar\tens\theta-\extd X\bar\tens \extd X),\quad \nabla\theta=0,
\quad\nabla\extd X=\theta\tens\extd X-\extd X\tens\theta.\]
This is manifestly not metric compatible with $\sigma=\id$ and indeed the connection does not have a classical origin.
\end{example}

\subsection{Metric compatible connection on a Cayley graph}

Here we focus on the canonical Euclidean metric in Section~3.3 but more general $\nabla$, seeking an actually metric compatible connection for it. 

Since the constraint in Lemma~3.8 involves squares in the graph, we look at Cayley graphs where there are many of these. Let $G$ be a finite group and let $\CC\subseteq G\setminus\{e\}$ be a set of generators closed under group inversion (here $e$ is the group identity element). In this case $\CC$ defines a left-translation invariant differential structure $\Omega^1$ on the finite group defined as in Section~3.1 from the induced bidirectional, connected regular Cayley graph. Here vertices are $V=G$ and the edges are $E=\{x\to xa\ |\ a\in \CC\}$. In this setting one has left-invariant 1-forms in terms of which the Euclidean metric takes the form
\[ \cg=\sum_a e_a\bar\tens e_{a^{-1}},\quad \(e_a,e_b\)=\delta_{a^{-1},b};\quad e_a=\sum_{x\in G}\omega_{x\to xa}.\]
Here $\Omega^1$ is a free module over $k(G)$ with basis $\{e_a\}$ and is also right-translation invariant precisely when $\CC$ is ad-stable. Note also that $\theta=\sum_{a\in\CC} e_a$. 

\begin{proposition} Let $x\to y \to z$ be a 2-step in the Cayley graph on a finite group and suppose that the generating set $\CC$ is closed under inverses and self-conjugation. Then there is an induced 2-step $x\to x y^{-1}z\to z$ and we define
\[  \includegraphics[scale=.8]{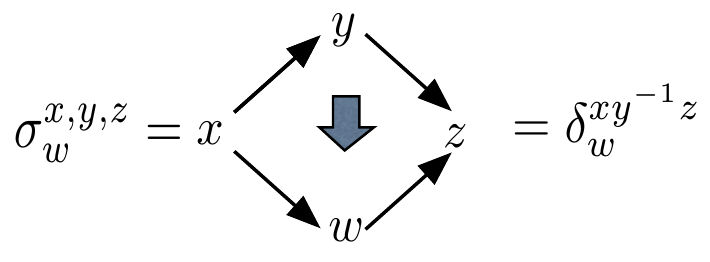}\]
Then $(\sigma,0)$ provides a metric compatible connection for the Euclidean metric,
\[ \nabla e_a=\sum_{b\in\CC}e_b\bar\tens(e_a-e_{b^{-1}ab}),\quad \sigma(e_a\bar\tens e_b)=e_b\bar\tens e_{b^{-1}ab}\]
where $\sigma$ is the standard crossed-module braiding on $\CC$ and obeys the braid relations. Moreover, the resulting  Laplace-Beltrami operator is $\Delta=2L$.
\end{proposition}
\proof Let $y=xa$ and $z=yb$. Define $w=x y^{-1}z$ then $w=xb$ and $z=xab=w(b^{-1}ab)$ so defines another 2-step $x\to w\to z$. We use this to define $\sigma$ and then verify the condition on $\sigma$ in Lemma~3.8. This can be done from the diagrammatic form or more conventionally
 \[ \forall \quad  \begin{array}{rcl}x&\longrightarrow&z\\ \downarrow & &\uparrow\\ v&\to&w \end{array},\quad \sum_{y:x\from y\to w}\sigma^{y,x,z}_w\sigma^{x,y,w}_v= \sum_{y:x\from y\to w}\delta^{yx^{-1}z}_w\delta^{xy^{-1}w}_v=\delta^z_v\]
 since the first delta function sets $y=wz^{-1}x$. From the first part of the proof we know that $w\to z\to x$ implies that $w\to y\to x$ is another 2-step hence a suitable $y$ exists in the sum. Since the graph is bidirectional we are not concerned about the reversal of some of the arrows in this explanation. Finally, we convert the results to the basis of left-invariant forms. Thus,
 \[ \sigma(\omega_{x\to y}\bar\tens\omega_{y\to z})=\omega_{x\to xy^{-1}z}\bar\tens\omega_{xy^{-1}z\to z}\]
 which in view of the above translates to $\sigma$ on the $\{e_a\}$ basis. This then gives $\nabla$ as $\nabla e_a=\theta\bar\tens e_a-\sigma(e_a\bar\tens\theta)$. Finally, the Laplace-Beltrami operator follows from Proposition~3.7 since $\gamma_{x,y}=2$ since only  $w=xy^{-1}z$ contributes in the sum and indeed $x\to w$ by our observation at the start of the proof.  \endproof
 
 We call this the `Maurer Cartan connection' for reasons that will emerge when we look at their geometry in the next section.

\begin{example} The bidirectional line graph with vertices $\Z$ is a Calyey graph with $\CC=\{\pm 1\}\subset \Z$. The basic forms $e_\pm=e_{\pm 1}$ are as in Example~3.9 but now, since the group is Abelian, $\nabla e_\pm=0$ and $\sigma={\rm flip}$ on the $e_\pm$. Here,  for example $\sigma(\omega_{x\to x+1}\bar\tens\omega_{x+1\to x+2})=\omega_{x\to x+1}\bar\tens\omega_{x+1\to x+2}$ but  $\sigma(\omega_{x\to x+1}\bar\tens\omega_{x+1\to x})=\omega_{x\to x-1}\bar\tens\omega_{x-1\to x}$. In geometrical terms $\nabla\extd X=\nabla\theta=0$ and $\sigma={\rm flip}$ on $\extd X,\theta$. 

\end{example}

\begin{example} The graph below is a Calyey graph in two different ways, namely $\CC=\{u,v,w\}\subset S_3$ where $u=(12),v=(23),w=(13)$, and  $\CC=\{1,3,5\}\subset \Z_6$, as indicated by the different labellings,
\[ \includegraphics[scale=.5]{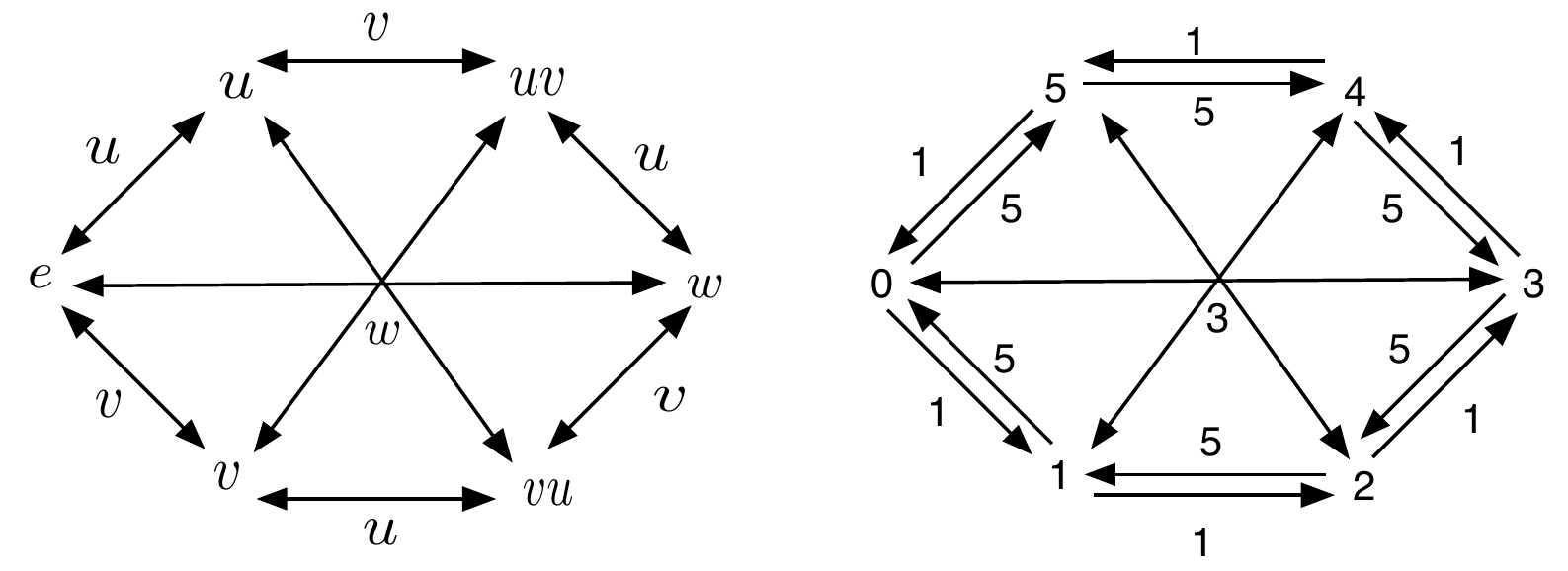} \]
These group structures define different connections $(\sigma, 0)$ compatible for the same metric. Thus for the Cayley graph on the left, 
\[ \nabla e_u=e_v\bar\tens (e_u-e_w)+e_w\bar\tens(e_u-e_v),\quad\sigma(e_u\bar\tens e_v)=e_v\bar\tens e_w\]
and similarly by cyclic permutation of $u,v,w$. For the Cayley graph on the right,
\[ \nabla e_a=0,\quad \sigma(e_a\bar\tens e_b)=e_b\bar\tens e_a\]
for all $a,b\in\Z_6$. Note that other connections are certainly possible in the 2nd case as there are squares not of the form in Proposition~3.10. 
\end{example} 

\begin{example} We analyse all connections compatible with the Euclidean metric in the case of a single square. There is only one nontrivial square up to rotation or reflection, hence for all $x\to y\to z$,
\[ \sigma^{x,y,z}_w=s_{x,y,z}\delta^y_w+t_{x,y,z}\delta_w^*\]
 for some coefficients $s_{x,y,z},t_{x,y,z}\in k$. Here $*$ is the unique vertex connecting to $x$ other than $y$. The different cases of the quadratic equation for $\sigma$ in Lemma~3.8  become
\[ t_{v,x,z}s_{x,v,w}+t_{z,x,z}t_{x,z,w}=0,\quad t_{v,x,z}t_{x,v,w}+t_{z,x,z}s_{x,z,w}=1,\quad \forall \begin{array}{rcl}x&\to& z\\ \downarrow & &\uparrow\\ v&\to & w\end{array}\]
where the vertices are distinct, and 
\[ s_{v,x,z}s_{x,v,x}+s_{z,x,z}t_{x,z,x}=0,\quad s_{v,x,z}t_{x,v,x}+s_{z,x,z}s_{x,z,x}=1,\quad \forall \begin{array}{rcl}x&\to& z\\ \downarrow & &\uparrow\\ v&\to & x\end{array}\]
where $v\ne z$. In both cases the sum in the quadratic equation is over $y=v,z$. We simplify these equations by looking at the `symmetric' case
\[ s^{x,y,x}=s,\quad t^{x,y,x}=t,\quad s^{x,y,z}=S,\quad t^{x,y,z}=T\]
for any  $x\to y$ or any $x\to y\to z$ with $z\ne x$, respectively. Then our equations reduce to
\[ TS+tT=0,\quad T^2+tS=1,\quad Ss+st=0,\quad St+s^2=1.\]
We have two 1-parameter classes of solutions
\begin{enumerate}
\item $ s\ne 0,\quad s^2-t^2=1,\quad T^2=s^2,\quad S=-t$ 
\[ {\rm e.g.}\quad s=T=1,\quad t=S=0\]
matches Proposition~3.10 for the Cayley graph of $\{(0,1),(1,0)\}\subset \Z_2\times\Z_2$
\item   $s=0,\quad tS=1$ 
\[ {\rm e.g.}\quad s=T=0,\quad t=S=1\]
matches Proposition~3.10 for the Cayley graph of $\{1,3\}\subset\Z_4$.\end{enumerate}
\end{example}

As this example shows,  solutions provided by presentation as a Cayley graph are special within the class of metric compatible bimodule connection. 

\subsection{Connections of permutation type}

Inspired by the special connections of the preceding sections, we propose a class of connections of combinatorial interest and study these further. Both the metric compatible  connection we found on Cayley graphs and the canonical connection on any graph are examples.

\begin{definition} We say that a bimodule connection  given by $(\sigma,\alpha)$ is of {\em permutation type} if for all $x,z$ the matrix $\sigma^{x,\cdot,z}_\cdot$ is a permutation.
\end{definition}
In this case our diagrammatic representation of $\sigma$ becomes an operation on the space of 2-arcs from fixed endpoints. Thus for fixed $x,z$ the operation sends $x\to y\to z$ to $x\to\sigma_{x,z}(y)\to z$ for some permutation $\sigma_{x,z}$. Here $\sigma_{x,z}(y)$ is the unique $w$ that makes $\sigma^{x,y,z}_w$ non-zero.

It is immediate from Proposition~3.7 that a bimodule connection of permutation type obeys $(\ ,\ )\sigma=(\ ,\ )$ {\em iff} 
\[ g_{\sigma_{x,x}(y)\to x}=g_{y\to x},\quad\forall y\to x.\] For example, this holds and also $\Delta=2L$ by Proposition~3.7 for every permutation connection and the Euclidean metric. Likewise, it is immediate from Lemma~3.8 that a permutation type connection defined by  $(\sigma,0)$ is metric compatible {\rm iff}
\[  \forall \quad  \begin{array}{rcl}x&\longrightarrow&z\\ \downarrow & &\uparrow\\ v&\to&w \end{array},\quad \left (\sigma_{\sigma_{x,w}^{-1}(v),z}(x)=w \ \Leftrightarrow v=z\right),\quad g_{z\to w}=g_{x\to \sigma^{-1}_{x,w}(z)}.\]

Moreover, 

\begin{proposition} If a bimodule connection on a graph is of permutation type:  \begin{enumerate}
\item $\sigma(\theta\bar\tens\theta)=\theta\bar\tens\theta$ and there is a canonical prolongation $\Omega^2$ by Lemma~2.2.
\item  In this case the  connection corresponding to $(\sigma,0)$ has 
\[ T_\nabla\omega_{x\to y}=-\sum_{z:y\to z}(\omega_{x\to y}\wedge\omega_{y\to z}+\omega_{x\to \sigma_{x,z}(y)}\wedge\omega_{ \sigma_{x,z}(y)\to z})\]
\[ R_\nabla\omega_{x\to y}=-\sum_{z,w:y\to z\to w}\omega_{x\to \sigma_{x,z}(y)}\wedge\omega_{\sigma_{x,z}(y)\to \sigma_{\sigma_{x,z}(y),w}(z)}\bar\tens\omega_{\sigma_{\sigma_{x,z}(y),w}(z)\to w}\]\end{enumerate}
\end{proposition}
\proof (1) From the definition of $\sigma$, 
\[ \sigma(\theta\bar\tens\theta)=\sum_{\begin{array}{rcl}x\to & y&\to z\\ \searrow & w & \nearrow\end{array}}\sigma^{x,y,z}_w\omega_{x\to w}\tens \omega_{w\to z}=\sum_{x\to y\to z}\omega_{x\to \sigma_{x,z}(y)}\tens \omega_{\sigma_{x,z}(y)\to z}=\theta\bar\tens\theta\]
since in doing the sum we may replace the sum over $y$ by a sum over $\sigma_{x,z}(y)$ for each $x,z$ that are 2 steps apart. We then construct the 2nd order differential structure by Lemma~2.2. (2) We compute $T_\nabla$ and $R_\nabla$ from Theorem~2.1.
\endproof

\begin{corollary} For a Cayley graph on a finite group with generating set $\CC$ stable under conjugation and inversion,   the `Maurer-Cartan connection' in Proposition~3.10 induces a canonical prolongation to $\Omega^2$ and obeys
\[  T_\nabla e_a=-e_a\wedge\theta-\sum_b e_b\wedge e_{b^{-1}ab},\quad  R_\nabla e_a=0\]
\end{corollary}
\proof Clearly the $\sigma$ in Proposition~3.10 is of permutation type since for each $x,z$ fixed the map $y\mapsto xy^{-1}z$ is a bijection. Hence Theorem~2.1 and Lemma~2.2 apply and the first of these provides the formula for the torsion. The curvature vanishes by Lemma~2.2 since $\sigma$ obeys the braid relations as already remarked. \endproof

Note that the torsion vanishes when the group is Abelian. This reproduces general results for finite group geometry and its Maurer-Cartan bimodule connection in \cite{BegMa:star} but now within a more general graph theory approach  independent of quantum group methods used there. The terminology is justified because on a Lie group there is a canonical `Maurer-Cartan connection' typically with torsion but zero curvature.  Note that \cite{BegMa:star} also finds, by Mathematica, a metric compatible {\em and} torsion free bimodule connection on $S_3$ while other works \cite{Ma:rief,NML} provide cotorsion and torsion free connections on $S_3,A_4$, all of these with curvature, whereas Proposition~3.15 suggests that the correct role of permutation-type connections is  to define $\Omega^2$ (and higher) as typified but not limited to the Cayley graph example. We can then look for a different metric-compatible torsion free or `Levi-Civita' connection or cotorsion and torsion free `generalised Levi-Civita' connection with respect to that as the further noncommutative Riemannian geometry. 

For completeness, we nevertheless compute the `nonstandard Ricci' curvature of Section~2.5 for permutation-type connections $(\sigma,0)$. This comes out for the Euclidean metric as
\[ S_{\nabla,\cg}=-\kern-20pt \sum_{x,z,w:\, y\in {\rm Fix}(\sigma_{x,z})\atop x\to y\to z\to w}\omega_{x\to \sigma_{y,w}(z)}\bar\tens\omega_{\sigma_{y,w}(z)\to w}\]
For the Euclidean metric compatible `Maurer-Cartan' connection coming from a Cayley graph in Proposition~3.10  and Corollary~3.16, we have
\[ S_{\nabla,\cg}=-\theta\bar\tens\theta.\]
This is non-zero because of the non-classical lift $\tilde R_\nabla$ which projects to $R_\nabla=0$  but which is not itself zero. The `nonstandard scalar curvature' 
\[ (\ ,\ )S_{\nabla,\cg}=-(\theta,\theta)=-|\CC|\] 
has constant value in this Maurer-Cartan case. 

\section{2nd order extension of inner calculi and the edge Laplacian}

In this section we extend the standard or `classical' differential structure on graphs in Section~3.1 by the graph Laplacian in a manner similar to the noncommutative differential structure on any Riemannian manifold\cite{Ma:bh}. Here both the initial `classical' differential structure and the extended one are noncommutative. In the process we naturally extend the Laplacian to 1-forms on the original graph. Our constructions are general but we then specialise to the graph case.

\subsection{Laplacians on inner differential algebras} 

We suppose that $\Omega^1$ is a differential structure on an algebra $A$ equipped with 2nd order operator  $\Delta: A\to A$. Recall from Section~2.3 that this comes along with a  bimodule map $\<\ ,\ \>:\Omega^1\bar\tens\Omega^1\to A$ as the associated `bivector field'.  Here $\Delta$ may not be a Laplace -Beltrami operator associated to any connection and $\<\ ,\  \>$ may not be related to the inverse of any metric. 

\begin{proposition} Given a second order operator $(\Delta,\<\ ,\ \>)$ on a differential algebra  $(A,\Omega^1,\extd)$, we let $\tilde\Omega^1=\Omega^1\oplus A\theta'$. This forms a  possibly non-surjective differential structure with bimodule structure and exterior derivative
\[[\theta',f]=0,\quad f\bullet \omega= f\omega,\quad \omega\bullet f=\omega f+\lambda \<\omega,\extd f\>\theta',\quad  \tilde\extd f=\extd f + {\lambda\over 2}(\Delta f)\theta'\]
for all $f\in A$, $\omega\in\Omega^1$. Here $\lambda\in k$ is a parameter.
\end{proposition}
\proof The proof is identical to that of \cite[Lemma~2.1]{Ma:bh} as nowhere were used commutativity of $A$ or of the original differential structure on it. There is, however, a change of notation. \endproof

Note that we have not required here that $f\tens g\mapsto f\tilde\extd g$ is surjective, which depends on $\Delta$ (one can require this).  However,  the extended differential structure is connected  if the original one is.  Subtracting the extended bimodule structures,
\[ [\omega,f]_\bullet=[\omega,f]_{\Omega^1}+\lambda\<\omega,\extd f\>\theta'\]
for the commutator in the new differential structure in terms of the original.  We denote by $\widetilde{\tens}$ the tensor product over $A$ defined with the new bimodule structures. If $\nabla$ is a connection on $\Omega^1$ we define $\nabla_\omega=(\<\omega,\ \>\tens\id)\nabla$ for all $\omega\in\Omega^1$. 

\begin{proposition} Let  $\nabla$ be a left connection on $\Omega^1$.
 \begin{enumerate}\item There is a left module map 
 \[ \phi:\Omega^1\bar\tens\Omega^1\to \tilde\Omega^1\widetilde{\tens}\tilde\Omega^1,\quad \phi(\omega\bar\tens\eta)=\omega\widetilde{\tens}\eta-\lambda\theta'\widetilde{\tens}\nabla_\omega\eta,\quad\forall\omega,\eta\in\Omega^1\]
\item If $\Delta$ extends to $\Omega^1$ in such a way that
\[ \Delta(f\omega)=(\Delta f)\omega+f\Delta\omega+2\nabla_{\extd f}\omega,\quad \forall f\in A,\ \omega\in\Omega^1\]
then
\[ \tilde\nabla\omega=\phi(\nabla\omega)+{\lambda\over 2}\theta'\widetilde{\tens}(\Delta-K)\omega,\quad\forall\omega\in\Omega^1\]
is a left connection for any left-module map $K$.
\end{enumerate}\end{proposition}
\proof The proof is identical to that in \cite[Lemma~2.2]{Ma:bh} and \cite[Lemma~2.3]{Ma:bh} respectively as nowhere were  commutativity of $A$ or the original differential structure, nor more than a left-connection and the stated property used. \endproof

This is as far as we will take the  theory generalising \cite{Ma:bh} for the moment. We turn to the other direction, starting with an inner calculus.

\begin{proposition} Suppose that $\Omega^1$ on $A$ is inner with generator $\theta$ and $\<\  ,\ \>$ is a bimodule map.
\begin{enumerate}
\item  $\Delta f=2 \<\theta,\extd f\>,\quad\forall f\in A$ fulfils  the conditions needed in Proposition~4.1.
\item The extended differential structure is inner with the same generator $\theta$.
\item  If $\sigma,(\ ,\ )$  are bimodule maps on $\Omega^1\bar\tens\Omega^1$  such that 
\[[ (\ ,\ )\sigma(\theta\bar\tens\theta),f]=0,\quad\forall f\in A,\quad \<\ ,\ \>={1\over 2}\(\ ,\ \)(\id+\sigma)\]
then $\(\ ,\ \)\nabla\extd=\Delta$ in (1), where  $\nabla$ is the connection corresponding to $(\sigma, 0)$ in Theorem~2.1.
\item In this case $\Delta=(\ ,\ )\nabla^2$ extends $\Delta$ to 1-forms in the manner needed in Proposition~4.2.
\item Explicitly, for $\omega\in\Omega^1$,
\[ \Delta\omega=((\ ,\  )\tens\id)[\theta\bar\tens\theta\bar\tens \omega-(\sigma_{23}+\sigma_{12}\sigma_{23})(\theta\bar\tens\omega\bar\tens\theta)+\sigma_{12}\sigma_{23}\sigma_{12}(\omega\bar\tens\theta\bar\tens\theta)]\]
\item If $\sigma(\theta\bar\tens\theta)=\theta\bar\tens\theta$ then $\Delta(f\theta)=(\Delta f)\theta$ for all $f\in A$ and $\Delta\theta=0$.
 \end{enumerate}
\end{proposition}
\proof (1) We verify, working in the original differential structure
\[ \<\theta,\extd(fg)\>=\<\theta,f\extd g\>+\<\theta,(\extd f)g\>=\<\theta f,\extd g\>+(\Delta f)g=\<\extd f,\extd g\>+f\Delta g+(\Delta f)g\]
for all $f,g\in A$. (2) We verfy the commutator in the extended differential structure 
\[ [\theta,f]=[\theta,f]_{\Omega^1}+\lambda\<\theta,\extd f\>\theta'=\extd f+{\lambda\over 2}(\Delta f)\theta'=\tilde\extd f\]
for all $f\in A$ using part (1). (3) If we are given $\sigma,(\ ,\ )$ as stated then 
\begin{eqnarray*} (\ ,\ )\nabla\extd f&=&(\ ,\ )(\theta\tens\extd f-\sigma(\extd f\tens\theta))=2\<\theta,\extd f\>-(\ ,\ )\sigma(\extd f\bar\tens\theta+\theta\bar\tens\extd f)\\  &=&2\<\theta,\extd f\>+(\ ,\ )\sigma(f\theta\bar\tens\theta- \theta\bar\tens\theta f)=2\<\theta,\extd f\>\end{eqnarray*}
for all $f\in A$. (4) Immediate from (\ref{lapder1}).  (5) We compute
\begin{eqnarray*}\Delta\omega&=&\(\ ,\ )_{12}( \nabla\tens\id+(\sigma\tens\id)(\id\tens\nabla))(\theta\bar\tens\omega-\sigma(\omega\tens\theta)\\
&=&(\ ,\ )_{12}[\theta\bar\tens\theta\bar\tens\omega-\sigma(\theta\bar\tens\theta)\bar\tens\omega-\theta\bar\tens\sigma(\omega\bar\tens\theta)+\sigma(\sigma_1\bar\tens\theta)\bar\tens\sigma_2\\
&&\quad\quad\quad+\sigma_{12}(\theta\bar\tens\theta\bar\tens\omega-\theta\bar\tens\sigma(\omega\bar\tens\theta)-\sigma_1\bar\tens\theta\bar\tens\sigma_2+\sigma_1\bar\tens\sigma(\sigma_2\bar\tens\theta))]\end{eqnarray*}
where $\sigma_1\bar\tens\sigma_2=\sigma(\omega\bar\tens\theta)$ is notation. (6) Under the stated assumption we have $\nabla\theta=0$ in Theorem~2.1 and hence $\Delta\theta=0$ by (4). We then use the product rule also in (4). \endproof

The significance of part (1) is that the Laplacian has a deeper origin as the `partial derivative' adjoint to the 1-form $\theta$, a principle which is invisible in classical geometry since there can be no direction $\theta$. A special case is where $(\ ,\ )\sigma=(\ ,\ )$ and $(\theta,\theta)$ is central in which case the condition in part (3) holds with $\<\ ,\ \>=(\ ,\ )$.  We will only use this special case in what follows.

\subsection{Edge Laplacian on bidirected graphs} 

We now compute how the above looks for our differential structures on graphs, which we know to be inner.  We start at the level of Section~3.2  with general inverse metrics $(\ ,\ )$ as in Section~3.2 but under the constraints $(\ ,\ )\sigma=(\ ,\ )$ as in Proposition~3.7 for connections of the form $(\sigma,0)$, in both cases for simplicity. As in Section 3.2 onwards, we look only at bidirected graphs. 

\begin{proposition} The Laplace-Beltrami on $\Omega^1$ or `edge Laplacian' on a bidirected graph with connection $(\sigma,0)$ obeying $(\ ,\ )\sigma=(\ ,\ )$ is
\[ \Delta \omega_{x\to y}=\left(\sum_{z:z\to x}g_{z\to x}\right)\omega_{x\to y}-2\kern-15pt\sum_{w,z:\ \begin{array}{c}x\to y\\ \downarrow\quad\quad\downarrow\\ w\to z\end{array}}\kern -10pt g_{x\to w}\sigma^{x,y,z}_w\omega_{w\to z}+\kern-15pt\sum_{w,z,s: \begin{array}{rcl}&\nearrow w\searrow& \\ x&\to y\to&z\\ &\searrow s\swarrow&\end{array}}\kern -10pt g_{w\to x}\sigma^{x,y,z}_w\sigma^{w,z,s}_x\omega_{x\to s}\]
\end{proposition}
\proof We observe that 
\begin{equation}\label{thetaxy} (\theta,\omega_{x\to y})=g_{x\to y}\delta_y,\quad (\omega_{x\to y},\theta)=g_{y\to x}\delta_x\end{equation}
and  note that the Proposition~4.3 part (5) under our assumption simplifies to 
\begin{equation}\label{Deltaw} \Delta\omega=(\theta,\theta)\omega-2(\ ,\ )_{12}\sigma_{23}(\theta\bar\tens\omega\bar\tens\theta)+(\ ,\ )_{12}\sigma_{23}\sigma_{12}(\omega\bar\tens\theta\bar\tens\theta)\end{equation}
which we straightforwardly compute on $\omega_{x\to y}$. \endproof

Next we specialise  to bimodule connections $(\sigma, 0)$ of permutation type in Section~3.5. From  Proposition~3.15 we have $\sigma(\theta\bar\tens\theta)=\theta\bar\tens\theta$ and hence by Proposition~4.3 part (6) we have 
\[ \Delta\theta=0,\quad \Delta (f\theta)=(\Delta f)\theta\]
for all functions $f$ (so the spectrum includes the spectrum of $\Delta$ on functions). We also know from Section~3.5 that for permutation connections and the Euclidean metric the Laplacian is (twice) the usual graph Laplacian $L$. So for all permutation connections the edge Laplacian with Euclidean metric extends the usual graph Laplacian. Also, for the Euclidean metric the coefficient in first term in $\Delta$ in Proposition~4.4 is the degree function $\deg(x)$. The following example is in this class. 

\begin{corollary}  For a Cayley graph on a finite group with $\CC$ stable under self-conjugation, the  `Maurer-Cartan' connection  in Proposition~3.10 and the Euclidean metric give Laplace-Beltrami on left-invariant forms 
\[ \Delta e_a=2(|\CC| e_a-\sum_{b\in \CC} e_{b^{-1}ab}),\quad\forall a\in \CC.\]
\end{corollary}
\proof This is of permutation type so the preceding comments apply. The formula is more readily computed  from $\nabla e_a$ in Proposition~3.10 and (\ref{Deltaw}). \endproof

This is zero when the Cayley graph has an abelian group underlying it. Thus the Laplace-Beltrami on 1-forms can tell apart the two different Cayley structures on the graph in Example~3.12. For the $S_3$ Cayley graph the Laplacian on the vector space of  left-invariant 1-forms has eigenvalues $\{6,6,0\}$, over $\C$ i.e.  strictly positive aside from the mandatory zero mode $\theta$. 

\subsection{Canonical edge Laplacian on graphs}

For the remainder of the section we focus on the simplest case of the theory above, namely the canonical Euclidean metric and the canonical connection in Section~3.3. Here $\nabla$ is given in our classification by $(\sigma,0)$ where $\sigma=\id$, and is trivially of permutation type. Hence, as explained after Proposition~4.4, the Laplace-Beltrami on 1-forms extends the usual graph Laplacian on functions up to a normalisation. We call it the {\em canonical edge Laplacian} $\Delta$ of a graph. Clearly in this case,  from Proposition~4.3, it is given by
\[ \Delta\omega=(\theta,\theta)\omega-2(\theta,\omega)\theta+(\omega,\theta)\theta,\quad\forall \omega\in\Omega^1\]
and takes the explicit form
 \[ \Delta \omega_{x\to y}=\deg(x)\omega_{x\to y}-2\sum_{z: y\to z}\omega_{y\to z}+\sum_{z: x\to z}\omega_{x\to z}.\]
given by Proposition~4.4 setting $\sigma^{x,y,z}_w=\delta^y_w$ and $g_{x\to y}=1$ for all $x,y,z,w$. The edge vectors $\omega_{x\to y}$ form a basis and play a role like the $\delta_x$ for scalers. The vector space on which it acts is of dimension twice the number of undirected edges. 

\begin{corollary} In the $\{\omega_{x\to y}\}$ basis the column sums of the matrix of the canonical edge Laplacian are zero. The graph is regular iff the row sums are also zero. 
\end{corollary}
\proof The matrix $\CL$ here is labelled by the directed edges of the bidirected graph according to $\Delta\omega_{x\to y}=\sum_{w\to z}\CL_{x\to y,w\to z}\omega_{w\to z}$. The column sums as a row vector correspond to $\Delta\theta$ and hence vanish due to $\Delta\theta=0$. Meanwhile, the sum of the row corresponding to $x\to y$ is given by the sum of the coefficients of $\Delta\omega_{x\to y}$ which from the formula is $2(\deg(x)-\deg(y))$. Hence these all vanish iff  the graph is regular. \endproof

\begin{example} For the $m$-gon graph, $m>1$, we label the edges modulo $m$ calling those edges going clockwise $\omega^+_i$ and those going anticlockwise $\omega^-_i$ (say). The canonical edge Laplacian takes the form
\[\Delta\omega^\pm_i=3\omega^\pm_i-2\omega^\mp_i+\omega^\mp_{i\mp 1}-2\omega^\pm_{i\pm 1}.\]
Over $\C$, apart from the mandatory $0$ eigenvalue for $\theta$ the matrix of $\Delta$ has strictly positive eigenvalues in the range
$(0,8]$ (with maximal eigenvalue 8 precisely for $m$ even). Explicitly, there are $m$  eigenvectors $e^{2\pi\imath p j\over m}\theta$ with corresponding eigenvalue $8\sin({\pi p\over m})^2$ for $p=0,1,2,\cdots {m-1}$. Here $e^{2\pi\imath p j\over m}$  denotes  a function on the set of vertices labelled by $j=0,1,\cdots,m-1$ with $\omega^+_j=\omega_{j\to j+1}$ and $\omega^-_j=\omega_{j+1\to j}$. There are also $m$ eigenvectors 
\[ e^{2\pi \imath p j\over m}\sum_i (2\omega^+_{i}-\omega^+_{i-1}-2\omega^-_{i-2}+\omega^-_{i-1}),\quad p=0,1,2,\cdots,m-1\]
of eigenvalue $2$, but not necessarily independent from the previous set. Similarly over $\R$. The canonical edge Laplacian matrix for $m=3$ and  in basis order $\omega_i^+,\omega_i^-$ is
\[\left(
\begin{array}{cccccc}
 3 & -2 & 0 & -2 & 0 & 1 \\
 0 & 3 & -2 & 1 & -2 & 0 \\
 -2 & 0 & 3 & 0 & 1 & -2 \\
 -2 & 1 & 0 & 3 & 0 & -2 \\
 0 & -2 & 1 & -2 & 3 & 0 \\
 1 & 0 & -2 & 0 & -2 & 3
\end{array}
\right).\]
\end{example}

\begin{example} Over $\R$ the canonical edge Laplacian for another regular graph is:
\[ \includegraphics[scale=.3]{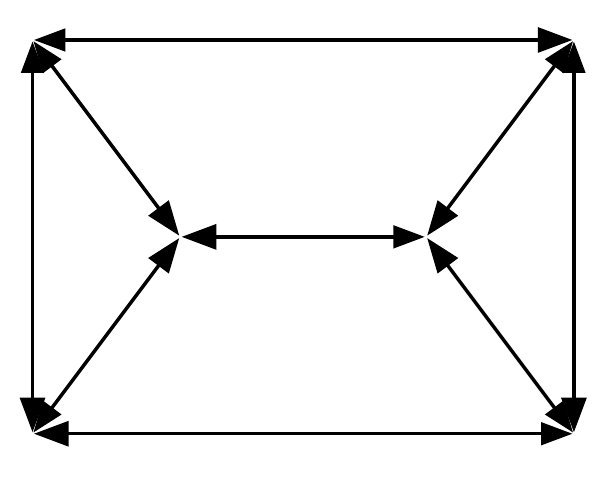}\quad {\rm Spec}(\Delta)=\{10(2),6(2),4, 3(12),0\}\]
where multiplicities are in brackets. We see that the eigenvalues are again strictly positive aside from the mandatory $0$-eigenvector $\theta$. Similarly for a non-regular graph:
\[ \includegraphics[scale=.3]{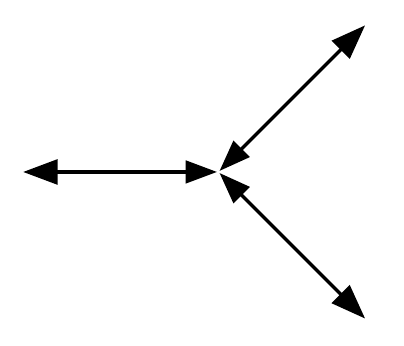}\quad {\rm Spec}(\Delta)=\{8,3(2),2(2),0\}\]
Note that this is quite different from the spectrum $\{6(2),0\}$ of $\Delta$ on functions (twice the usual  graph Laplacian on vertices)  applied to  the line graph, in this case a triangle. 
\end{example}

We now show that the canonical edge Laplacian matrix, although not necessarily symmetric, is nevertheless always positive semidefinite in the sense of its eigenvalues. In fact we determine the spectrum.

\begin{theorem} On a connected graph with vertex set $V$:

(1) The canonical edge Laplacian $\Delta$ has spectrum
\[{\rm Spec}(\Delta)=2\, {\rm Spec}(L)\cup  \{\deg(x)\ ({\rm with\ algebraic\ multiplicity\ }\deg(x)-1)\ |\ x\in V\} \]
where $L$ is the usual graph Laplacian on vertices and $\deg(x)$ is the degree of $x$.  

(2) Over $\R$, the canonical edge Laplacian is triangularisable  and its spectrum is strictly positive aside from one zero mode $\theta$.

(3) Over $\R$, the canonical edge Laplacian is fully diagonalisable  if the two stated parts of the spectrum  are disjoint. 
\end{theorem}
\proof (1) For each vertex $x\in V$ we let 
\[ \theta_x=\sum_{y:x\to y}\omega_{x\to y}\in \Omega^1\]
and in these terms we have
\[ \Delta\omega_{x\to y}=\deg(x)\omega_{x\to y}-2\theta_y+\theta_x.\]
Next, we find
\[ \Delta\theta_x=\sum_{y:x\to y}\Delta\omega_{x\to y}=\sum_{y:x\to y}(\deg(x)\omega_{x\to y}-2\theta_y+\theta_x)=2 (\deg(x)\theta_x-\sum_{y:x\to y}\theta_y)\]
The right hand side here is exactly the same matrix as that of $2L$ but now on the basis $\theta_x=\delta_x\theta$. We now construct a basis with sub-basis $\{\theta_x\ |\ x\in V\}$. Here  $\theta_x\in \span\{\omega_{x\to y}\}=\delta_x\Omega^1$ and we let $H_x$ be a chosen complement such that $\delta_x\Omega^1=k\theta_x\oplus H_x$. These blocks $H_x$ have dimension $\deg(x)-1$ and to be concrete we choose basis $\{\omega_{x\to y}\ |\ y\ne y_0\}$ where we leave out one $\omega_{x\to y_0}$. This requires a choice of one distinguished edge coming out of each vertex. Our basis consists of these bases of each $H_x$ and then the sub-basis $\{\theta_x\}$. Then we see that $\Delta$ has the form
\[ \left(\begin{array}{cccccccccc|cccccccccc}
                                           &     d(x_1)1_{d(x_1)-1}        &   &     & \vdots &        &  & & & &    &              &  &\vdots& &           &   &           & \cr         &       &   \ddots        &     & \vdots &        &  & & & &    &              &  &\vdots& &           &   &           & \cr

                                          0&\cdots &0 & d()&0     &\cdots&0&0&\cdots&0& 0 & \cdots & 0 &1    &0&\cdots &-2&\cdots&0\cr
                                           0&\cdots &0 & 0     &d()&\cdots&0&0&\cdots&0& 0 & \cdots & 0 &1     &0&\cdots &0&\cdots&-2\cr
                                            &              &   &     & \vdots &  \ddots      &  & & & &    &              &  &\vdots& &           &   &           & \cr
                                           0&\cdots  &0 & 0 & \cdots&   & d()&0&\cdots&0& 0  & \cdots &0 &1    & -2& \cdots&0&       0 & 0\cr
                                                 &              &   &     & \vdots &        &  & & \ddots& &    &              &  &\vdots& &           &   &           & \cr

                                            &              &   &     & \vdots &        &  & & &   d(x_m)1_{d(x_m)-1}  &    &              &  &\vdots& &           &   &           &\cr \hline
                                             &              &   &     & 0 &        &  & & & &    &              &  &2 L & &           &   &           & \cr
                                        
                                           \end{array}\right)\]
for an enumeration of the vertices and with $d=\deg$ for brevity. We show the rows for the generic block $H_x$ in detail. The -2's appear one on each row and in different columns corresponding the the different $y$ such that $x\to y$ and $y\ne y_0$. From this  triangular form the algebraic spectrum is immediate (from the characteristic polynomial). Although we have been concrete, there are other natural choices of $H_x$ that one could use, with a similarly block-triangular form of $\Delta$. 

(2) Over $\R$, it is known\cite{Chung} that $L$ is already fully diagonalisable with just one zero mode in the case of a connected graph. This immediately carries forward to $\Delta$ since  $\deg(x)\ge 1$ for a connected graph and also means in view of the above that $\Delta$ is triangularisable.

(3)  In our case $L$ acts on the space spanned by the $\{\theta_x\}$ and as mentioned already has a full set of eigenvectors. As $\Delta$ restricts to $2L$ we obtain  $m=|V|$ eigenvectors of $\Delta$. Next consider extending a basis element $\omega_{x\to y}$ of $H_x$ to an eigenvector $\omega_{x\to y}+v$ of $\Delta$ with eigenvalue $\deg(x)$, where $v$ is in the space spanned by the $\{\theta_z\}$. We need to solve $ \Delta(\omega_{x\to y}+v)=\deg(x)(\omega_{x\to y}+v)$ but the left hand side is $\deg(x)\omega_{x\to y}-2\theta_y+\theta_x+2Lv$ hence we need to solve
\[ (2L-\deg(x))v=2\theta_y-\theta_x.\]
This will have a solution for all $y\ne y_0$ if the two parts of the spectrum in (1) are disjoint. At least in this case we will obtain a full set of eigenvectors and $\Delta$ will be fully diagonalisable. \endproof

The Examples~4.8 are both `generic' in the sense of the two parts of the spectrum being disjoint, and hence fully diagonalisable. On the other hand, it is not the case that the canonical edge Laplacian is always fully diagonalisable and the best we can do in general is to put it in triangular form, as the following example shows.

\begin{example} An example of a graph with non-diagonalisable canonical edge Laplacian is
\[ \includegraphics[scale=.3]{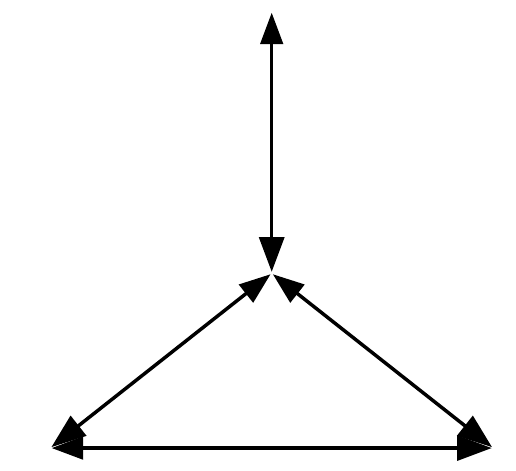} \quad {\rm Spec}(\Delta)=\{8,6,3(2),2(3),0\}\]
Here the 2nd part of the spectrum is $\{3(2),2(2)\}$ from the degrees of the vertices, hence $2L$ has spectrum $\{8,6,2,0\}$ so the two parts the spectrum are not disjoint. Correspondingly, $\Delta$  has only two independent eigenvectors of eigenvalue 2.
\end{example}

The $m$-gon in Example~4.7 is mixed: the two parts of the spectrum are disjoint precisely when $m\neq 0$ mod 6 and this case is diagonalisable. One can check that when $m=6,12,18$ (for example) the canonical edge Laplacian is not diagonalisable. 
Note that Theorem~4.9 says that disjointness of the two parts of the spectrum is sufficient for diagonalisability. We have not proven that this is necessary although examples such as the above would be consistent with this.

\begin{lemma} On a Cayley graph with finite group and generating set $\CC$ the canonical edge Laplacian  has eigenvectors of the form\begin{enumerate}
\item   $f\theta$ with eigenvalue $\mu$,  $\forall f$ such that  $\Delta f=\mu f$.
\item $\sum v^ae_a$ with eigenvalue $|\CC|$, $\forall \{v^a\}$ such that $ \sum_{}(2R_{a^{-1}}-1) v^a= 0$. Here $(R_af)(x)=f(xa)$ and the sums are over $a\in\CC$.\end{enumerate}
\end{lemma}
\proof  The first type of eigenvector is a restatement of $\Delta(f\theta)=(\Delta f)\theta$. Recall that $\Delta f=2L f$ on functions. Next, from the abstract formula for the canonical edge Laplacian and $(\theta,e_a)=(e_a,\theta)=1$ we see immediately that
\[ \Delta e_a=|\CC|e_a-\theta.\]
Then
\begin{eqnarray*} \Delta(fe_a)&=&(\Delta f)e_a+f\Delta e_a+2\nabla_{\extd f}e_a\\
&=&(\Delta f+n f)e_a-f\theta+2(\extd f,\theta)e_a-2(\extd f,e_a)\theta=
|\CC| fe_a -(2R_{a^{-1}}f-f)\theta\end{eqnarray*}
using the results of Proposition~4.2.  This gives the second type of eigenvector.
 \endproof

This includes  Example~4.7 and the other circulants or product of circulants, which can now be understood using characters. Such graphs are those that can be expressed as induced from Cayley graphs on finite Abelian groups. Thus, given a circulant graph of $m$ vertices, label them by $0,1,\cdots, m-1\in\Z_m$ and let $\CC$ be the column labels where the top row of the adjacency matrix is 1. The adjacency matrix at entry $(i,j)$ is then 1 iff $j-i\in \CC$ and as the matrix is symmetric we see that $\CC$ is closed under negation. Hence the graph is induced by a Cayley graph structure on $\Z_m$; the converse is equally clear. The following then applies, which also explicitly constructs the eigenvalues and eigenvectors where possible.

\begin{proposition} On a Cayley graph with finite Abelian group $G$ and generating set  $\CC$  (i.e. a connected  circulant or product of such) the canonical  edge Laplacian over $\C$ has $|G|$ eignevalues and eigenvectors of the form (1) in Lemma~4.11 coming from the usual graph Laplacian on $G$ and $(|\CC|-1)|G|$ eigenvectors of the form (2) in Lemma~4.11 with eigenvalue $|\CC|$.  The two collections are not necessarily linearly independent. \end{proposition}
\proof  The set $\CC\subseteq G\setminus\{e\}$ is required to generate so that the induced graph is connected and to be closed under inverses so that the graph is bidirected.  The first type of eigenvalues and eigenvectors in Lemma~4.11 are merely inherited from the  Laplacian on functions (twice the usual graph Laplacian) and for completeness we begin by diagonalising that (this is not new but we do it in our notations). Let $\hat G$ be the group of unitary 1-dimensional representations and $\chi\in \hat G$. Then  $\chi(xa)=\chi(x)\chi(a)$ for all $x\in G$  hence 
\[ \Delta\chi=2L\chi=2\sum_{a\in\CC}(\chi-\chi(a)\chi)=2(\sum_{a\in\CC}(1-\chi(a)))\chi.\]
Now, if $a^2=e$ then $\chi(a)^2=1$ and $\chi(a)=\pm 1$. Hence terms of this form contribute 4 or 0 to the eigenvalue. If $a\ne a^{-1}$ then these both occur in the sum and together they contribute
\[ 2(2-\chi(a)-\chi(a^{-1}))=2(2-\chi(a)-\chi(a)^{-1})=8({\chi(a)^{1\over 2}-\chi(a)^{-{1\over 2}}\over 2\imath})^2\ge 0\]
since the character values are necessarily roots of unity. Again we have  0  iff $\chi(a)=1$.  In the entire sum then we obtain 0 iff $\chi(a)=1$ for all $a\in\CC$. But in this case $\chi(g)=1$ for all $g\in G$ since $\CC$ generates i.e. only for the trivial representation. Otherwise we confirm that the eigenvalues are positive (as known in general for any connected graph). As $|\hat G|=|G|$ this fully diagonalises $\Delta$ on functions  (with eigenvalues algebraic integers since they are integer sums of roots of unity) and by Lemma~4.11 it also gives us $|G|$ eigenvectors $\chi\theta$ of the edge Laplacian as we run over $\chi\in \hat G$. The new part is the second type of eigenvector; we consider $v^a=\chi\mu^a$ where $\mu^a\in\C$ are constants and we compute
\[ \sum_{a\in \CC}(2R_{a^{-1}}-1)(\chi\mu^a)=\chi\sum_{a\in\CC}(2\chi(a^{-1})-1)\mu^a\]
which vanishes precisely when $\mu=(\mu^a)$ is perpendicular to the vector   in $\C^{|\CC|}$ with entries $2\chi(a^{-1})-1$ in a basis labelled by $a\in\CC$. Hence there are $|\CC|-1$ independent directions for $\mu$ for each $\chi\in\hat G$, hence $(|\CC|-1)|G|$ eigenvectors in total of the second type in Lemma~4.11, all with eigenvalue $|\CC|$. In total we have constructed  $|\CC||G|$ eigenvectors, i.e. $2|{\rm undirected\ edges}|$  by the handshaking lemma. This is the number of directed edges and hence we have potentially diagonalised the edge Laplacian, i.e. if they are all linearly independent. The case of the $m$-gon where $m$ is a multiple of 6 shows that they need not be. \endproof

The Laplacian on functions is similarly diagonalised by matrix elements of irreducible representations when $G$ is nonAbelian   and $\CC$ is ad-stable, while the edge Laplacian is more complicated. We also remark that for any regular finite graph we can similarly colour the edges in each direction. We have seen in Proposition~3.2 that following the edges of each colour (they will form non-intersecting loops) provides a basis $\{\omega_a\}$. These generalise the bases $\{e_a\}$ on a Cayley graph and it is shown in \cite{MaRai} that such a colouring has the interpretation of a vielbein of `permutation type'. More of the geometry of such colourings is a direction for study elsewhere.

It remains further to study and apply the edge Laplacian. It also remains to find reasonable conditions for the `extended connection' $\tilde\nabla$ in Proposition~4.2 to be a bimodule connection. It would also be of interest to connect to Hodge theory as in the example \cite{MaRai2} and to the Killing form and exterior algebra on finite groups\cite{Ma:per,LMR}.



\begin{thebibliography}{99}

\bibitem{BegMa:star}
E.J. Beggs and S. Majid, *-Compatible connections in noncommutative Riemannian geometry,  J. Geom. Phys. 25 (2011) 95-124

\bibitem{BrzMa} T. Brzezinski and S. Majid, Quantum differentials and the q-monopole revisited, Acta Appl. Math. 54 (1998) 185-232


\bibitem{Chung} F. Chung, {\em Spectral Graph Theory}, CBMS Regional Conference Series in Mathematics, No. 92, AMS 1997

\bibitem{Connes} A. Connes.
\newblock {\it Noncommutative geometry}, Academic Press, 1994

\bibitem{DVM1}
 M. Dubois-Violette and T. Masson, On the first-order operators in bimodules, Lett. Math. Phys. 37 (1996) 467Ð474
 
 \bibitem{DVM2}M. Dubois-Violette and P.W. Michor, Connections on central bimodules in noncommutative differential geometry, J. Geom. Phys. 20, 218 Ð232, 1996.
 
\bibitem{Hall} P. Hall, On Representatives of Subsets, {\em J. London Math. Soc.} 10 (1935) 26Ð30

\bibitem{HelNes} P. Hell and J. Nesetril, {\em Graphs and Homomorphisms}, Oxford University Press, 2004

\bibitem{LMR}
J. Lopez Pena, S. Majid and K. Rietsch, Lie theory of finite simple groups and the generalised Roth conjecture, 38pp., arXiv:1003.5611 (math.QA)

\bibitem{Ma:book}
S. Majid, {\em Foundations of Quantum Group Theory}, Cambridge University Press, 1995

\bibitem{Ma:rief}
S.~Majid.
\newblock Riemannian geometry of quantum groups and finite groups with
  nonuniversal differentials.
\newblock {\em Commun. Math. Phys.}, 225:131-170, 2002.

\bibitem{MaRai}
S. Majid and E. Raineri. 
\newblock Moduli of quantum Riemannian geometries on $\le$ 4 points, J. Math. Phys. 45 (2004) 4596-4627.


\bibitem{MaRai2}
S. Majid and E. Raineri, Electromagnetism and gauge theory on the permutation group $S_3$, J. Geom. Phys. 44 (2002) 129-155

\bibitem{Ma:per}
S. Majid, Noncommutative differentials and Yang-Mills on permutation groups $S_N$, Lect. Notes Pure Appl. Maths 239 (2004) 189-214, Marcel Dekker

\bibitem{Ma:bh}
S. Majid, Almost commutative Riemannian geometry, I: wave operators, 39pp, Commun. Math. Phys, (2011), in press.

\bibitem{Mou}
J. Mourad, Linear connections in noncommutative geometry, Class. Quantum Grav. 12 (1995) 965 Ð 974.

\bibitem{NML}
F. Ngakeu, S. Majid and D. Lambert, Noncommutative Riemannian Geometry of the Alternating Group $A_4$, J. Geom. Phys. 42 (2002) 259-282

 
\end{thebibliography}
\end{document}